\documentclass[12pt]{article}
\usepackage{latexsym}
\usepackage{graphicx}
\usepackage{verbatim}
\usepackage{amssymb}
\usepackage{amsmath}
\title{Generalized dimensions of images of measures under Gaussian processes}
\author{Kenneth Falconer \\
\small{{\it Mathematical Institute,
University of St~Andrews, North Haugh, St~Andrews,}} \\
\small{{\it Fife, KY16~9SS, Scotland }}\\
\and
Yimin Xiao\thanks{Research of Y. Xiao was partially supported by
NSF grants DMS-1006903 and DMS-1309856.}\\
\small{{\it Department of Statistics and Probability,
Michigan State University,}}\\
\small{{\it East Lansing, MI 48824, USA}}}

\def\bbbr{{\mathbb R}}
\def\R{{\mathbb R}}

\def\bbbn{{\mathbb N}}
\def\bbbz{{\mathbb Z}}
\def\a{\alpha}
\def\la{\lambda}

\def\de{\delta}
\def\ga{\gamma}

\def\eps{\varepsilon}
\def\si{\sigma}

\newtheorem{theo}{Theorem}
\newtheorem{prop}[theo]{Proposition}

\newtheorem{lem}[theo]{Lemma}
\newtheorem{cor}[theo]{Corollary}

\newcommand\pkd{\mbox{\rm dim}_{\rm P}\,} 
\newcommand\hdd{\mbox{\rm dim}_{\rm H}\,} 
\newcommand\qd{D_{q}\,}
\newcommand\lqd{\underline{D}_{q}\,}
\newcommand\uqd{\overline{D}_{q}\,}

\newcommand{\bi}{x} 
\newcommand{\bj}{y} 
\newcommand{\bu}{u} 
\newcommand{\bv}{v} 
\newcommand{\bw}{w} 
\newcommand{\E}{{\sf E}} 
\renewcommand{\P}{{\sf P}} 
\renewcommand{\O}{{\cal O}} 
\newcommand{\Orb}{{\rm Orb}} 
\newcommand{\be}{\begin{equation}} 
\newcommand{\ee}{\end{equation}} 

\allowdisplaybreaks

\begin{document}
\maketitle

\begin{abstract}
\noindent We show that for certain Gaussian random processes and fields $X:\bbbr^N  \to \bbbr^d$,
$$
\qd (  \mu_X)  = \min\Big\{d, \, \frac 1{\alpha}\qd (\mu)\Big\} \quad \mbox{ a.s.,}
$$
for an index $\alpha$ which depends on H\"{o}lder properties and strong local nondeterminism
of $X$, where $q>1$, where $\qd$ denotes generalized $q$-dimension  and where $\mu_X$ is the image of the measure
 $\mu$ under $X$. In particular this holds for index-$\alpha$ fractional Brownian
motion, for fractional Riesz-Bessel motions and for certain infinity scale fractional Brownian motions.
\end{abstract}

\section{Introduction}
\setcounter{equation}{0}
Dimensions of images of sets under stochastic processes have been studied for many years.
The Hausdorff dimension of the image or sample path of Brownian motion $X:\bbbr_+ \to \bbbr^d$
is almost surely equal to
$$ \hdd X(\bbbr_+) = \min \big\{ d, 2\big\},$$
where $\hdd$ denotes Hausdorff dimension, see L\'{e}vy  \cite{Le}, with the exact gauge function for the dimension
 established by Ciesielski and Taylor \cite{C-Tay62} for $d \ge 3$ and by Ray \cite{Ray63} and
Taylor \cite{Tay} for $d=2$. Similar questions were subsequently studied for other processes,
notably for sample paths of stable L\'evy processes, see \cite{Tay67}, and for fractional Brownian motion,
see \cite{Ad1,Ad,Kah,ManVN,Tal95,Tal98}. There are several comprehensive surveys of this work \cite{Ad,Kah,Tay2,Xi04}
which contain many further references.

A more general, but very natural, question is to find the almost sure dimensions of the image
$X(E)$ of a Borel set $E \subseteq \bbbr^N$ under a process $X:\bbbr^N  \to \bbbr^d$, in
terms of the dimension of $E$. In particular, Kahane \cite{Kah} showed that
\be
 \hdd X(E) = \min \Big\{ d, \frac{\hdd E}{\alpha}\Big\} \quad \mbox{ a.s.} \label{dimset}
 \ee
if $X$ is index-$\alpha$ fractional Brownian motion (which reduces to standard Brownian motion when
$\alpha = \frac{1}{2}$ and $N=1$).

The corresponding question for packing dimension $\pkd$\!, where dimensions of images of sets can
behave in a more subtle manner, was not answered until rather later, when Xiao \cite{Xi1}
showed that for index-$\alpha$ fractional Brownian motion,
$$ \pkd X(E) = \frac{\dim_P^{\alpha d} E} {\alpha} \quad \mbox{ a.s.},$$
where $\dim_P^{s} E$ is the `packing dimension profile' of $E$, a notion introduced in
connection with linear projections of sets by Falconer and Howroyd \cite{FH}, and which
is defined in terms of a certain $s$-dimensional kernel.

In recent years, many other dimensional properties of the range, graph, level sets and
images of given sets have been studied for a wide range of random processes, see
\cite{Ad,Kho,Xi97,Xi07,Xi09b} for surveys of this work. 

It is natural to study dimensional properties of images of measures under random processes
or fields in an analogous way to images of sets. For $\mu$ a Borel measure on $\bbbr^N$  
and  $X:\bbbr^N  \to \bbbr^d$, the random image measure $ \mu_X$ on $\bbbr^d$ is defined by
$$ 
\mu_X(A) = \mu\{x: X(x) \in A \}, \quad A \subseteq \bbbr^d.
$$
When $\mu$ is the Lebesgue measure on $\bbbr^N$ and $X$ is a Gaussian process, the properties 
of the corresponding image measure $\mu_X$ have played important roles in studying the exact 
Hausdorff measure functions for the range, graph and level sets of $X$ \cite{Tal95,Xi97}.
For more general Borel measures $\mu$, one can look at the almost sure Hausdorff and packing dimensions of the measures
(given by the minimal dimension of any set with complement of zero measure);  indeed, by supporting suitable
measures on sets, this approach is often implicit in the set dimension results mentioned above.
Explicit results for Hausdorff and packing dimensions of image measures under a wide range
of processes are given in \cite{SX}, with dimension profiles again key in the packing dimension cases.

However, the singularity structure of a measure may be very rich, and multifractal analysis
in various forms has evolved to exhibit this structure as a function or spectrum; for general
discussions see, for example, \cite{Fa,Har,Lau,P}.  In this paper we consider generalized
$q$-dimensions which reflect the asymptotic behaviour as $r \searrow 0$ of the $q$th-moment sums
$M_{r}(q) =  \sum_{C} \mu (C)^{q}$  over the mesh cubes $C$ of side $r$ in $\bbbr^N$. It
will be convenient for us to work with the equivalent $q$th-moment integrals $
\int\mu(B(x,r))^{q-1}d\mu(x)$, where $B(x,r)$ is the ball centre $x$ and radius $r$,
see Section 2 or \cite{Lau,PS} for further details of $q$-dimensions.

Our main results are natural measure analogues of (\ref{dimset}) for the generalized $q$-dimension.
Thus, for example, for Gaussian processes $X:\bbbr^N  \to \bbbr^d$ which are strongly locally
$\alpha$-nondeterministic and which satisfy an $\alpha$-H\"{o}lder condition, and for a
compactly supported probability measure $\mu$ on $\bbbr^N$, we show that
$$
\lqd (  \mu_X)  = \min\Big\{d, \, \frac 1{\alpha}\lqd (\mu)\Big\} \quad \mbox{ a.s.},
$$
where $\lqd$ denotes (lower) generalized $q$-dimension. Such processes include index-$\alpha$
fractional Brownian motion, fractional Riesz-Bessel motions and infinity scale fractional Brownian motion.
We restrict attention to `larger' moments, that is where $q>1$.

Typically, upper bounds for the generalized dimensions follow easily from the almost sure
H\"{o}lder continuity of the sample paths.  Lower bounds are more elusive because $X:\bbbr^N
\to \bbbr^d$ may map many disparate points into a small ball. Simplifying things somewhat, 
we need to obtain upper bounds for  $\mu_X(B(X(y),r))^{q-1}$ by bounding their expectations. Taking $n$ to be the integer with
$n\leq q<n+1$, this requires bounding the probability that the images under $X$ of sets of points
$\{x_1,\ldots,x_n,y\}$ in $\bbbr^N$ are within distance $r$ of each other, that is, essentially,  
\be
\P\big\{
 |X(y)-X(x_i)|  \leq r \mbox{ for all } 1\leq i \leq n\big\}. \label{probballs}
 \ee
For this we use the strong local non-determinism of  $X$, roughly that the conditional variance 
 $\mbox{var}\big(X(x_i)| X(x_j)\  (j \neq i)\big)$ is comparable 
 to the unconditional variance 
 $\mbox{var}\big(X(x_i)-X(x_j)\big)$ where $x_j$ is the point (with  $j\neq i$) closest to $x_i$.
This enables us to obtain a bound for (\ref{probballs}) in terms of an expression $\phi(x_1,\ldots,x_n,y)$  
which reflects the mutual distance distribution of the points $\{x_1,\ldots,x_n,y\}$.
To estimate this, we integrate $\phi$ with respect to $\mu$ for each $x_1,\ldots,x_n$, then 
a further integration over $y$ bounds the generalized $q$-dimension in terms of integrals of $\phi$.

This reduces the problem to  estimating  `multipotential' integrals of the form
$$\int\bigg[
 \int \cdots \int
 \phi(x_1,\ldots,x_n,y) d\mu(x_1)\ldots d\mu(x_n)\bigg]^{(q-1)/n}d\mu(y). $$
 Two stages are then needed to estimate such integrals. In Section 3 we introduce a device that allows us to replace integrals over Euclidean space by more tractable ones over an ultrametric space, so that the integral  becomes an infinite sum over the vertices of an $(n+1)$-ary tree. We then estimate this in Section 4 using an induction process over `join' vertices of tree, at each step using H\"{o}lder's inequality to estimate the relevant sums over increasingly large sets of vertices on the tree. Obtaining these estimates turns out to be significantly more awkward when $q$ is non-integral and we need to associate a particular path in the tree with the fractional part of $q$.

\section{Main definitions and results}
\setcounter{equation}{0}
\setcounter{theo}{0}

This section details generalized $q$-di\-men\-sions of 
measures and the random processes that we will be concerned with, to enable us to state our main results.

We first review generalized $q$-di\-men\-sions of measures, which are the
basis for the `coarse' approach to multifractal analysis, that is the approach
based on `box sums', as opposed to the `fine' approach based on Hausdorff or packing
dimensions, see \cite{Fa,Gr,Har,Lau,P}
for various treatments.
The {\it $r$-mesh cubes} in $\bbbr^N$ are the cubes in the family
$${\cal M}_{r} = \big\{ [i_{1}r,(i_{1}+1)r) \times \cdots \times
[i_{N}r,(i_{N}+1)r) \subseteq \bbbr^N: i_{1}, \ldots , i_{N} \in \bbbz\big\}.$$
Throughout the paper $\mu$ will be a finite Borel measure of bounded
support on $\bbbr^{N}$.
For $q >0$ and $r>0$  the $q$th-{\it moment sums} of $\mu$ are given by
$$
M_{r}(q) = \sum_{C\in  {\cal M}_{r}} \mu (C)^{q},
$$
For $q>0, q\neq 1$, we define the {\it lower} and {\it upper generalized
q-dimensions} or {\it R\'{e}nyi dimensions} of $\mu$ to be
\be
\lqd (\mu) = \liminf_{r \rightarrow 0}\frac{\log  M_{r}(q)}{(q-1) \log r}
 \quad \mbox{and} \quad
\uqd (\mu) = \limsup_{r \rightarrow 0}\frac{\log  M_{r}(q)}{(q-1) \log
r}.\label{3.b}
\ee
If, as happens for many measures, $\lqd (\mu)=\uqd (\mu)$, we write $D_{q}(\mu)$ for
the common value which we refer to as the {\it generalized q-dimension}. 
When we just write $\qd (\mu)$ it is implicit that the generalized $q$-dimension exists.
Note that the definitions of $q$-dimensions are independent of the origin
and coordinate orientation chosen for the mesh cubes.

There are useful integral forms of
$\lqd$ and $\uqd$.  For $q>0, q\neq 1$,
\begin{eqnarray}
\lqd (\mu) &=& \liminf_{r \rightarrow 0}\frac{\log \int
\mu(B(x,r))^{q-1}d\mu(x)}{(q-1) \log r} \label{intdef1}\label{deflgd}\\
\mbox{and} \quad
\uqd (\mu) &=& \limsup_{r \rightarrow 0}\frac{\log \int
\mu(B(x,r))^{q-1}d\mu(x)}{(q-1) \log r};\label{intdef2}
\end{eqnarray}
see \cite{Lau} for $q>1$ and \cite{PS} for $0<q<1$.

It is easily verified that $\lqd (\mu)$ and $\uqd (\mu)$ are each
nonincreasing in $q$ and continuous (for  $q\neq 1$), and that
$0 \leq \lqd (\mu) \leq \uqd (\mu) \leq N$ for all $q$.

Let $\mu$ be a Borel probability measure
on $\bbbr^N$ with compact support, and let $f:\bbbr^N \to \bbbr^d$ be Borel measurable and bounded on the support of $\mu$.
 The {\it image measure} $ \mu_f$ of $\mu$ under
$f$ is the Borel measure of bounded support  defined by
$$
 \mu_f(A) = \mu\{x: f(x) \in A\}, \quad A \subseteq \bbbr^d.
$$
In particular,  for any measurable $g: \bbbr^d \to \bbbr_+$
$$
\int_{\bbbr^d} g(y) d \mu_f(y) = \int_{\bbbr^N} g(f(x) )d \mu(x).
$$

It follows easily from (\ref{intdef1}) and  (\ref{intdef2}) that if $f:\bbbr^N \to \bbbr^d$ is
$\alpha$-H\"{o}lder on compact intervals in $\bbbr^N $, where $0 < \alpha \le 1$, and $q>0, q \neq 1$, then
\be
\lqd ( \mu_f)   \leq \min\Big\{d,\, \frac 1 \alpha\, \lqd (\mu) \Big\} \quad \mbox{ and } \quad
\uqd ( \mu_f)   \leq \min\Big\{d,\, \frac 1 \alpha\, \uqd (\mu) \Big\}. \label{Holder}
\ee
These inqualities may be regarded as measure analogues of the relationships between
the Hausdorff, packing and  box dimensions of sets and their images under H\"{o}lder
mappings, see \cite{Fa}.

Now let $X:\bbbr^N \to \bbbr^d$ be a continuous random process or random field on
a probability space $(\Omega, {\cal F}, \P)$ and let $\E$ denote expectation.
Then for $\mu$ a Borel probability measure on $\bbbr^N$ with compact support, the random
image measure $ \mu_X$ on $\bbbr^d$  is defined by
$$
\mu_X(A) = \mu\{x: X(x) \in A\},\quad  A \mbox{ a Borel subset of } \bbbr^d.
$$
The main aim of the paper is to relate the $q$-dimensions of
$\mu_X$ and $\mu$ for suitable Gaussian processes.

Immediately from  (\ref{Holder}), if $X: \bbbr^N \to \bbbr^d$ almost surely has
a H\"{o}lder exponent $\alpha$ on a compact interval $K \subseteq \bbbr^N$ which
contains the support of $\mu$, then for $q>0, q \neq 1$
\be
\lqd ( \mu_X)   \leq \min\Big\{d,\, \frac 1 \alpha\, \lqd (\mu) \Big\} \quad \mbox{ and } \quad
\uqd ( \mu_X)   \leq \min\Big\{d,\, \frac 1 \alpha\, \uqd (\mu) \Big\} \mbox{ a.s. }  \label{holderbounds}
\ee
For many processes, Kolmogorov's continuity theorem, see for example \cite{RW},
provides a suitable H\"{o}lder exponent. Much of our effort will be devoted to obtaining inequalities in the opposite direction to those of (\ref{holderbounds}).

Henceforth we assume that $X:\bbbr^N \to \bbbr^d$ is a Gaussian random field defined by
\begin{equation}\label{def:X}
X(x) = \big(X_1(x), \ldots, X_d(x)\big),\qquad x \in \R^N,
\end{equation}
where $X_1, \ldots, X_d$ are independent copies of a mean zero Gaussian process $X_0:
\bbbr^N \to \bbbr$ with $X_0(0) = 0$ a.s.
We also assume that $X_0$  satisfies the following
Condition {\bf (C)}:
\begin{itemize}
\item[{\bf (C)}] For some $\delta_0>0$, let $\psi: [0, \delta_0) \to [0, \infty)$ be a
non-decreasing, right continuous
function with $\psi(0) =0$ such that, for some constant $C_{1}>0$,
$$
\frac{\psi(2r)} {\psi(r)} \le C_{1} \quad \mbox{ for all }\ r \in (0,\de_0/2).
$$
We assume:
\begin{itemize}
\item[{\bf (C1)}]\ There is a constant  $C_{2}>0$ such that for all $x,h \in \R^N$ with $|h| \le \de_0$,
$$
\E\bigl[\bigl( X_0(x + h)- X_0(x)\bigr)^2\bigr] \le C_{2} \psi(|h|).
$$
\item[{\bf (C2)}]\ For all $T > 0$, the process $X_0$ is {\it strongly locally
$\psi$-nondeterministic} on $[-T, T]^N$, that is,  there exist
positive constants $C_{3}$ and $r_0 $ such that for all $x \in [-T, T]^N$
and all $0 < r \le \min \{ |x|,\ r_0 \}$,
\begin{equation}\label{Eq:Cuzick}
{\rm Var} \bigl( X_0(x) \big | X_0(y) \, :\,  y \in [-T, T]^N,\ r \le |x - y| \le r_0
\bigr) \ge C_{3}\, \psi(r).
\end{equation}
\end{itemize}
\end{itemize}

The concept of local nondeterminism was first introduced by Berman \cite{Be} for
Gaussian processes and was subsequently extended by Pitt \cite{Pitt} to random fields.
The above definition of strong local $\psi$-nondeterminism is essentially due
to Cuzick and DuPreez \cite{CD} (who considered the case $N=1$).
For brief historical details
and various applications of strong local nondeterminism see \cite{Xi06,Xi07}.

Dimensional
properties of  Gaussian fields have been studied in \cite{Ad,Kho,Xi97,Xi07} in increasing
generality. The almost sure Hausdorff dimensions of certain random sets associated with a
random process $X$ given by (\ref{def:X}), such as the range $X([0, 1]^N)$, graph
Gr$X([0, 1]^N) = \{(x, X(x)): x \in [0, 1]^N\}$
and level sets $X^{-1}(z)= \{ x \in [0, 1]^N: X(x) = z\}$,
depend on the {\it upper index of} $\psi$ {\it at} 0 defined by
\begin{equation} \label{Def:phiup}
\a^* = \inf \Big\{\beta \ge 0: \lim_{r \to 0} \frac{\psi(r)}
{r^{2\beta}} = \infty\Big\}
\end{equation}
(with the convention $\inf \emptyset = \infty$). Analogously, we can define the
{\it lower index of } $\psi$  {\it at}  0 by
\begin{equation} \label{Def:philow}
\a_* = \sup \Big\{\beta \ge 0: \lim_{r \to 0} \frac{\psi(r)}
{r^{2\beta}}= 0\Big\}.
\end{equation}
Clearly, $0 \le \a_* \le \a^* \le \infty$. Under Condition {\bf (C)}, the indices $\a^*$ and $\a_*$ are uniquely 
determined by $X_0$. Xiao \cite{Xi09b} showed that the packing dimension of the range
$X([0, 1]^N)$ is determined by the lower index $\a_*$. Moreover, if $\psi(r)$ is regularly varying
at $r=0$ of index $2\alpha\in (0, 2]$, then $\a_* = \a^* = \alpha$.

We may now state the main result of the paper.
\begin{theo}\label{main}
Let $X: \R^N \to \bbbr^d$ be a Gaussian random field defined by $(\ref{def:X})$ and assume
that the associated random process $X_0$ satisfies Condition {\bf (C)} for some $\psi$. Let $q > 1$ and
let $\mu$ be a Borel probability measure on $\R^N$ with compact support.
\begin{itemize}
\item[(i)]\, If $ 0<\a_* = \a^*=\alpha < 1$, then
\begin{equation}\label{Eq:main1}
\lqd (\mu_X)  = \min\Big\{d, \, \frac 1{\alpha}\lqd (\mu)\Big\} \quad \mbox{a.s.}
\end{equation}
\item[(ii)]\, If the generalized $q$-dimension $D_q (\mu)$ of $\mu$ exists and $0 < \a_* \leq \a^* < 1$, then
\begin{equation}\label{Eq:main2}
\min\Big\{d, \, \frac 1{\alpha^*}D_q(\mu) \Big\} =  \lqd (\mu_X) \leq \uqd (\mu_X)
\leq\min\Big\{d, \, \frac 1{\alpha_*}D_q(\mu) \Big\}
\quad \hbox{a.s.}
\end{equation}
\item[(iii)]\, If the generalized $q$-dimension $D_q (\mu)$ of $\mu$  exists and $ 0<\a_*
= \a^*=\alpha < 1$, then the generalized $q$-dimension $D_q (\mu_X)$ of $\mu_X$ exists almost surely, and
$$\qd (\mu_X)  = \min\Big\{d, \, \frac 1{\alpha}\qd (\mu)\Big\} \quad \mbox{a.s.}$$
\end{itemize}
\end{theo}

There are several examples to which Theorem \ref{main} is readily applicable. The most important 
Gaussian random fields which satisfy Condition {\bf (C)}
are the fractional Brownian motions. Recall that, for $0<\alpha<1$, a real-valued {\it index-$\alpha$
fractional Brownian motion} $B^\alpha:  \R^N \to \bbbr$  is the centered Gaussian
random field with covariance function
\be
\E\big(B^\alpha(x) B^\alpha(y) \big) = \frac1 2\big( |x|^{2\alpha} + |y|^{2 \alpha}
- |x - y|^{2 \alpha}\big), \quad x,y \in \bbbr^N, \label{fbm}
\ee
introduced by Mandelbrot and Van Ness \cite{ManVN} for $N=1$. When $N > 1$ and $\alpha = 1/2$,
then $B^\alpha$ is L\'evy's $N$-parameter Brownian motion, see \cite[Chapter 18]{Kah}.
It follows that $\E \big[\big(B^\alpha(x)- B^\alpha(y) \big)^2\big] =
|x - y|^{2 \alpha}$ so $B^\alpha$ has stationary, isotropic increments
and is $\alpha$-self-similar. Strong local $\psi$-nondeterminism of $B^\alpha$, with
$\psi(r) = r^{2\alpha}$, follows from Lemma 7.1 of Pitt \cite{Pitt}, whose proof relies
on the self-similarity of $B^\alpha$. A different proof using Fourier analysis
can be found in \cite{Xi07}.

\begin{cor}[Fractional Brownian motion]\label{Coro:FBM}
Let $X: \R^N \to \bbbr^d$ be index-$\alpha$ fractional Brownian motion $(\ref{fbm})$. Let
$\mu$ be a Borel probability measure on $\bbbr^N$ with compact support.  Then for all $q>1$,
$$
\lqd (  \mu_X)  = \min\Big\{d, \, \frac 1{\alpha}\lqd (\mu)\Big\} \quad
\mbox{ and } \quad \qd (  \mu_X)  = \min\Big\{d, \, \frac 1{\alpha}\qd (\mu)\Big\}
\quad \mbox{a.s.}
$$
\end{cor}

Another example is {\it fractional Riesz-Bessel motion} $Y^{\ga, \beta}:
\R^N\to \bbbr$, with indices $\ga$ and $\beta$, introduced by Anh, Angulo and
Ruiz-Medina \cite{Anh99}. This is a centered Gaussian random field
with stationary increments and spectral density
\be
f_{\gamma, \beta}(\la) = \frac {c(\ga, \beta, N)} {|\la|^{2 \ga}
(1 + |\la|^2)^\beta}, \label{rb}
\ee
where $\ga$ and $\beta$ are constants satisfying
$$\textstyle  \beta + \ga - \frac{N} 2 >0, \quad 0 < \ga < 1 + \frac{N} 2$$
and $c(\ga, \beta, N)>0$ is a normalizing constant. These  random fields are important
for modelling as they exhibit long range dependence and intermittency simultaneously,
see \cite{Anh99}.

It may be shown that if $\gamma + \beta - \frac N 2 > 1$ then the sample function
$Y^{\ga, \beta}(x)$ has continuous (first order)  partial derivatives almost surely, see, for example, \cite{XX11}.
So the generalized dimensions of $\mu_X$ are the same as those of $\mu$. For
$0< \gamma + \beta - \frac N 2 < 1$, it is proved in \cite{Xi07} that $Y^{\ga, \beta}$
satisfies Condition {\bf (C)} with $\psi(r) = r^{2(\gamma + \beta - \frac N 2)}$, so applying 
Theorem \ref{main} leads to the following statement.

\begin{cor}[Fractional Riesz-Bessel motion]\label{Coro:rb}
Let $X: \R^N \to \bbbr^d$ be index-$(\gamma,\beta)$ fractional Riesz-Bessel motion with spectral density 
$(\ref{rb})$. Let  $\mu$ be a Borel probability measure on $\bbbr^N$ with compact support.  If
$\gamma+\beta -\frac{1}{2}N>1$ then $\lqd (  \mu_X) = \lqd (\mu)$ and  $\qd (  \mu_X)
= \qd (\mu)$ a.s. for all $q>0, q \neq 1$. If
$0 <\gamma+\beta -\frac{1}{2}N<1$ then
for all $q>1$,
$$
\lqd ( \mu_X)  = \min\Big\{d, \, \frac{ \lqd (\mu)}{(\gamma+\beta -\frac{1}{2}N)}\Big\}
\, \mbox{ and } \, \qd ( \mu_X)  = \min\Big\{d, \, \frac{ \qd (\mu)}
{(\gamma+\beta -\frac{1}{2}N)}\Big\}\quad a.s.
$$
\end{cor}

When $X_0$ has stationary and isotropic increments, $\a^{*}$ and $\a_*$
coincide with the upper and lower indices of $\si^2(h)$ defined in an analogous manner
to (\ref{Def:phiup})-(\ref{Def:philow}), where
$$
\si^2(h) = \E\bigl[\bigl( X_0(x + h)- X_0(x)\bigr)^2\bigr], \qquad
 x, h \in \R^N
$$
(by isotropy $\sigma^2(h)$ is a function of $\|h\|$). Many interesting examples
of Gaussian random fields with stationary increments which satisfy
condition  {\bf (C)} can be constructed, see \cite{EWX,LX11,Xi97,Xi07,Xi09}.

We recall a class of Gaussian random fields with $ \a_* < \a^*$, due to Clausel \cite{Clausel10}.
The approach is similar to the method for constructing L\'evy processes
with different upper and lower Blumenthal-Getoor indices \cite{BG61}.
We remark that, while Blumenthal and Getoor's indices are concerned with the asymptotic
behavior of $\si^2(h)$ as $\|h\|\to \infty$, we are interested in the behavior of
$\si^2(h)$ near $h=0$.

Let ${\bf H} = \{H_j, \, j \ge 0\}$ be a sequence of real numbers such that
\[
0 < \liminf_{j\to \infty} H_j \le \limsup_{j\to \infty}H_j < 1.
\]
A real-valued Gaussian random field $B_{\bf H}: \R^N \to \bbbr$
with stationary increments may be defined by the harmonizable representation:
$$
B_{{\bf H}} (t) = \sum_{j = 0}^\infty  \int_{D_j}
\frac{e^{i \langle t, \lambda \rangle}-1} {\|\lambda\|^{H_j + \frac N 2}}\,
dW( \lambda),
$$
where $D_0 = \{\lambda \in \R^N: \|\lambda \| < 1\}$ and $D_j =
\{\lambda \in \R^N: 2^{j -1}\le \|\lambda \| <  2^j\}$ for $j \ge 1$, see \cite{Clausel10}. Then
$B_{{\bf H}}$ is called the {\it infinity scale fractional Brownian motion}
with indices ${\bf H} = \{H_j, \, j \ge 0\}$. It is proved in \cite{EWX} that
$\alpha_* = \liminf\limits_{j\to \infty} H_j $ and, under an extra condition
on $\{H_j, \, j \ge 0\}$, we have $\alpha^* = \limsup\limits_{j\to \infty} H_j$.
To be more precise, let $\underline{H} = \liminf\limits_{j\to \infty} H_j$
and $\overline{H} = \limsup\limits_{j\to \infty} H_j$.
For each $\varepsilon \in ( 0, \underline{H})$, we define
a sequence $T_n =T_n(\varepsilon)$ as follows
\[
T_1 = \inf\{j: H_j \ge \overline{H}-\eps\},\quad \quad
T_2 = \inf\{j>T_1: H_j < \overline{H}-\eps\},
\]
and for all $k \ge 1$ define inductively
\[
T_{2k+1} = \inf\{j> T_{2k}: H_j \ge \overline{H}-\eps\}
\]
and
\[
T_{2k+2} = \inf\{j>T_{2k+1}: H_j < \overline{H}-\eps\}.
\]
If we assume that for every $\varepsilon \in ( 0, \underline{H})$,
\begin{equation}\label{Eq:Lacu2}
T_{2k+2} > \frac{(\overline{H}-\eps) ( 1- \underline{H} +\eps)}
{(\underline{H}-\eps) ( 1- \overline{H} +\eps)}\, T_{2k+1}
\end{equation}
for all $k$ large enough, then it can be verified that
$\alpha^* = \limsup\limits_{j\to \infty} H_j$,
see \cite{EWX}.

\begin{cor}[Infinity scale fractional Brownian motion]\label{Coro:ifsfbm}
Let $X: \R^N \to \bbbr^d$ be an infinity scale fractional Brownian motion with indices
${\bf H} = \{H_j, \, j \ge 0\}$, which satisfies $(\ref{Eq:Lacu2})$ for all $\varepsilon>0$
small enough. Let $q > 1$ and let $\mu$ be a Borel probability measure on
$\R^N$ with compact support.
\begin{itemize}
\item[(i)]\, If $ \lim_{n\to \infty} H_j = \alpha \in (0,  1)$, then
\begin{equation*}
\lqd (\mu_X)  = \min\Big\{d, \, \frac 1 {\alpha}\lqd (\mu) \Big\} \quad \hbox{a.s.}
\end{equation*}
\item[(ii)]\, If\, $0 < \liminf_{j\to \infty} H_j < \limsup_{j\to \infty}H_j < 1$,
then if $\mu$ has generalized $q$-dimension $D_q (\mu)$, we have
\begin{equation*}
\lqd (\mu_X)  = \min\Big\{d, \, \frac 1{\overline{H}}D_q(\mu) \Big\} \quad \hbox{a.s.},
\end{equation*}
where $\overline{H}= \limsup_{j\to \infty} H_j$.
\end{itemize}
\end{cor}

The proof of Theorem \ref{main} will be completed in Section 5. However, before 
we can make the necessary probability estimates we need some technical results, 
which are developed in Sections 3 and 4.

\section{Equivalent ultrametrics}
\setcounter{equation}{0}
\setcounter{theo}{0}

In this section we define, for each $n \in \bbbn$, an ultrametric $d$ on the unit cube
in $\bbbr^N$ and a finite number of translates $d_{a}$  (i.e. of the form $d_{a}(x,y) = d(x+a,y+a)$) 
and a constant $c\ge 1$, such that, given any set of $n$  points $x_1,\ldots, x_n \in [0,\frac{1}{2})^N$, 
there is a $d_{a}$ such that  $c^{-1} |x_i-x_j|  \leq d_{a}(x_i,x_j) \leq c|x_i-x_j|$ for all $1 \leq i,j \leq n$, 
where $|x-y|$ is the Euclidean distance between $x$ and $y$. 
Thus, on any given set of $n$ points,  one of the $d_{a}$  is equivalent
to the Euclidean metric in a uniform manner. We need this so we can replace
the Euclidean metric by an ultrametric when estimating the expectations that arise in
Section 5.

For $m \geq 2$ we construct a hierarchy of  $m$-ary subcubes of the unit cube $[0,1)^N$
in the usual way. For $k= 0,1,2,\ldots$ define the set of $k${\it -th level cubes}
\be
{\cal C}_k =\Big\{ [i_1 m^{-k}, (i_1 +1) m^{-k})\times\cdots\times   [i_N m^{-k}, (i_N +1) m^{-k})
: 0 \leq i_1,\ldots, i_N \leq m^k -1\Big\}.\label{cubes}
\ee
These cubes define an ultrametric $d$ on $[0,1)^N$ given by
$$
d(x,y) = m^{-k},  \qquad x,y \in [0,1)^N,
$$
where $k$ is the greatest integer such that $x$ and $y$ are in the same cube of $ {\cal C}_k$,
with $d(x,x)=0$.

Whilst it is easy to see that $|x-y| \leq \mbox{const} \,d(x,y)$, the opposite inequality is not uniformly valid. To
address this, we consider translates of $d$ to get a family of ultrametrics on $[0,\frac{1}
{2})^N$, from which we can always select one that will suit our needs.

Assume (to avoid the need for rounding fractions) that $m$ is even. Let ${\cal A}_m$ denote
the family of translation vectors:
\be
{\cal A}_m = \Big\{ \Big(\frac{j_1}{m-1}, \ldots, \frac{j_N}{m-1}\Big)
:  0 \leq j_1,\ldots,j_N \leq \frac{m}{2} -1   \Big\}. \label{translates}
\ee
For each $a \in {\cal A}_m $ define
\be
d_{a}(x,y) = d(x+a,y+a),  \qquad x,y \in [0,{\textstyle\frac{1}{2}})^N;\label{transum}
\ee
then $d_{a}$ is an ultrametric on $[0,\frac{1}{2})^N$ which we may think of as a translate of
$d$ by the vector $a$. (Note that the restriction on the indices $0 \leq j_l \leq \frac{m}{2} -1$
in (\ref{translates}) ensures that the $d_{a}$ are defined throughout $[0,\frac{1}{2})^N$.)

For $a \in {\cal A}_m$ we write ${\cal C}_k^a$ for the  cubes obtained by translating
the family ${\cal C}_k$ by a vector $-a$, that is
\be
{\cal C}_k^a = \big\{ C -a : C\in {\cal C}_k\} \label{cka};
\ee
thus $d_{a}(x,y)$ is also given by the greatest integer such $k$ such that $x$ and $y$ are
in the same cube of $ {\cal C}_k^a$.

\begin{prop}\label{propums}

\noindent (i) For all $a \in {\cal A}_m$ we have
$$
|x-y| \leq N^{1/2} d_a (x,y) \qquad x,y \in [0,{\textstyle\frac{1}{2}})^N.
$$

\noindent  (ii) Given $x,y \in [0,{\textstyle\frac{1}{2}})^N$, we have
\be
d_{a}(x,y) \leq 8m(m-1) |x-y| \label{ineq2}
\ee
for all except at most $N(\frac{m}{2})^{N-1} $ vectors  $a \in {\cal A}_m$.
\end{prop}

\noindent{\it Proof.}
(i) With $k$ the greatest integer such that $x+a$ and $y+a$ are in the same
$k$-th level cube of ${\cal C}_k$,
$$d_a (x,y) = d(x+a,y+a) = m^{-k} = N^{-1/2} N^{1/2} m^{-k} \geq N^{-1/2}|x-y|.$$

\noindent (ii)
We first prove (ii) in the case $N=1$.

Let $x \in  [0,\frac{1}{2})$ and let $k\geq 1$. We claim that, for all $j$ such
that $0 \leq j \leq  \frac{m}{2} -1 $ with at most one exception,
\be
\Big| x + \frac{j}{m-1} - i m^{-k} \Big| \geq \frac{1}{4m^k (m-1)} \quad
\mbox{ for all } i \in \bbbz.
\label{lbound}
\ee
Suppose, for a contradiction, that there are  $0 \leq j \neq j'\leq  \frac{m}{2} -1$
and $i,i' \in \bbbz$ such that both
$$ \Big| x + \frac{j}{m-1} - i m^{-k} \Big| <  \frac{1}{4m^k (m-1)} \quad \mbox{ and } \quad
\Big| x + \frac{j'}{m-1} - i' m^{-k} \Big| <  \frac{1}{4m^k (m-1)} .$$
Then
$$ \Big| \frac{j-j'}{m-1} - (i-i') m^{-k} \Big| <  \frac{1}{2m^k (m-1)} $$
so
$$
\big| (j-j')m^{k-1}(1+m^{-1} +m^{-2}+\cdots) - (i-i') \big |
= \Big| \frac{(j-j')m^k}{m-1} - (i-i') \Big| <  \frac{1}{2(m-1)}.
$$
Thus the integer $z = (j-j')m^{k-1}(1+m^{-1} \cdots +m^{-k+1}) - (i-i')$ satisfies
$$\Big| z+\frac{j-j'}{m-1}  \Big| =
\big|z+ (j-j')(m^{-1} +m^{-2}+\cdots) \big |  < \frac{1}{2(m-1)},
$$
which cannot hold for any $0 \leq j \neq j'\leq  \frac{m}{2} -1 $, proving the
claim (\ref{lbound}).

Now suppose  $x, y \in [0,\frac{1}{2})$ satisfy $|x-y| \leq 1/ (8m(m-1))$ and let
$k\geq 1$ be the integer such that
\be
\frac{1}{8m^{k+1}(m-1)} < |x-y| \leq \frac{1}{8m^{k}(m-1)}. \label{range}
\ee
For all $j$ such that (\ref{lbound}) holds for this $x$ and $k$, we have for all
$i \in \bbbz$\begin{eqnarray*}
\Big| y + \frac{j}{m-1} - i m^{-k}\Big|&\geq&\Big| x + \frac{j}{m-1} - i m^{-k}\Big| - |x-y| \\
&\geq & \frac{1}{4m^k (m-1)}  - |x-y| \\
&\geq &  |x-y| \qquad  \qquad  \qquad  \qquad    \mbox{(using (\ref{range}))}\\
&=& \Big| \Big(y + \frac{j}{m-1}\Big) -\Big(x + \frac{j}{m-1}\Big)\Big|.
\end{eqnarray*}
Hence, $y + j/(m-1)$ is in the same interval of ${\cal C}_k$ as $x + j/(m-1)$, so for
all $j$ such that (\ref{lbound}) holds, that is for all except at most one value of $j$,
$$ d_{j/(m-1)}(x,y)= d\big(x + j/(m-1),y + j/(m-1)\big)  \leq m^{-k} \leq 8m(m-1)|x-y|.$$
If $1/ (8m(m-1)) < |x-y| \leq \frac{1}{2}$ then $d_{j/(m-1)}(x,y) \leq 1 < 8m(m-1)|x-y|$
for all $j$, completing the proof of (ii) when $N=1$.

For $N\geq 2$, write, in coordinate form, $x=(x_{_1},\ldots,x_{_N}),y=(y_{_1},\ldots,y_{_N})
\in \bbbr^N$ and $a=(a_{_1},\ldots,a_{_N})\in {\cal A}_m$. Applying the result for $N=1$
to each coordinate,
$$ d_a(x,y) = \max_{1\leq l \leq N} d_{a_l}(x_{_l},y_{_l})
\leq 8m(m-1)\max_{1\leq l \leq N}|x_{_l}-y_{_l}|\leq 8m(m-1)|x-y|,$$
provided that $a_{_l}$ is not an exceptional value for the $1$-dimensional case for
any coordinate $l$. There are at most  $N(\frac{m}{2})^{N-1} $ such exceptional vectors
in $a\in {\cal A}_m$, otherwise (\ref{ineq2}) holds.
$\Box$

\medskip

\begin{cor}\label{corums}
Let $x_1,\ldots, x_n \in [0,\frac{1}{2})^N$. If $m>2n^2N$ then there exists $a \in {\cal A}_m$
such that
\be
N^{-1/2}|x_i -x_j| \leq  d_{a}(x_i,x_j) \leq  8m(m-1)|x_i-x_j| \label{pairs}
\ee
for all $1 \leq i,j \leq n$.
\end{cor}

\noindent{\it Proof.}
By Proposition \ref{propums}, for each pair $x_i,x_j $ there are at most $N(\frac{m}{2})^{N-1}$
vectors $a\in {\cal A}_m$ for which (\ref{pairs}) fails. Since there are a total of $(\frac{m}{2})^{N}$
vectors in  ${\cal A}_m$, this leaves at least  $(\frac{m}{2})^{N}- n^2 N(\frac{m}{2})^{N-1} $
vectors in  ${\cal A}_m$ such that  (\ref{pairs}) holds for all $1 \leq i,j \leq n$, and this number is positive if
$m>2n^2 N$.
$\Box$

\section{Integral estimates}\label{integralest}
\setcounter{equation}{0}
\setcounter{theo}{0}

This section is devoted to proving Theorem \ref{mainint} which bounds the multipotential integral (\ref{intest2}) which arises when estimating the $q$-dimensions of the image measures. A related procedure was used in
\cite{Fa5,Fa6} in connection with self-affine measures. We work with a {\it code space} or
{\it word space} on $M$ symbols, which we may identify with the vertices of the  $M$-ary rooted
tree in the usual way. In Section 5 this will in turn be identified in a natural way with the hierarchy of $m$-ary
cubes (\ref{cubes}) in the ultrametric construction of Section 3, with $M= m^N$.

For $k=0,1,2, \ldots$ let $I_{k}$ be the
set of all $k$-term sequences or words formed from the integers $1,2, \ldots, M$,
that is $I_{k} = \{ 1,2, \ldots, M \}^k$, with $I_{0}$ comprising the empty word $\emptyset$.
We write $|v| =k$ for
the length of a word $v \in I_k$.
Let
$I= \cup^{\infty}_{k=0} I_{k} $
denote the set of all finite words, and let $I_{\infty} = \{ 1,2, \ldots, M \}^\bbbn$ denote
the corresponding set of infinite words.
We write ${\bi}|_k$ for the {\it curtailment} of  $\bi\in I \cup I_\infty$ after $k$ terms,
that is the word comprising the initial $k$ terms of $\bi$.
For $v \in I$ and $\bi \in I \cup I_\infty$
we write $v \preceq \bi$ to mean that ${v}$ is an initial subword of
 ${\bi}$.  If ${\bi,\bj} \in I_{\infty} $  then
$\bi\wedge \bj$ is the maximal word such that both
$\bi\wedge \bj \preceq \bi$ and $\bi\wedge \bj \preceq \bj$.

We may topologise $I_{\infty}$ in the natural way by the metric
$d(\bi,\bj) = 2^{-|\bi \wedge \bj |}$ for distinct $\bi,\bj \in
I_{\infty}$ to make $I_{\infty}$  into a compact metric space, with
the {\it cylinders}
$C_v = \{\bj \in I_{\infty} : \bv \preceq \bj \}$ for  $v \in I$ forming a base
of open and closed neighborhoods of $I_{\infty}$.

 It is convenient to identify $I$ with the vertices of an $M$-ary
 rooted tree with root $\emptyset$. The edges of this tree join each vertex $v \in I$
 to its $M$  `children'  $v 1,\ldots,v M$.

The {\it join set} $\wedge(\bi_1,\ldots,\bi_n)$ of $\bi_1,\ldots,\bi_n \in I_\infty$
is the set of vertices $\{\bi_k \wedge\bi_{k'} : 1 \leq k \neq k' \leq n \}$. We say
that $\bu \in \wedge(\bi_1,\ldots,\bi_n)$ has  {\it multiplicity}  $r-1$ if $r$ is
the greatest integer for which there are distinct indices  $k_1,\ldots,k_r$ with
$\bi_{k_p}\wedge \bi_{k_q} = \bu$ for all $1 \leq p<q \leq r$, see Figure 1. (In the case of $M=2$
where $I$ is a binary tree, every vertex of a join set has multiplicity 1.) Counting
according to multiplicity, the join set $\wedge(\bi_1,\ldots,\bi_n)$ always comprises $n-1$
vertices of $I$. The {\it top vertex} $\wedge^T(\bi_1,\ldots,\bi_n)$ of a join set
is the vertex $\bv \in \wedge(\bi_1,\ldots,\bi_n)$ such that  $\bv \preceq \bi_l$ for
all $1 \leq l \leq n$.

\begin{figure}[t]
\begin{center}
\includegraphics[scale=0.42]{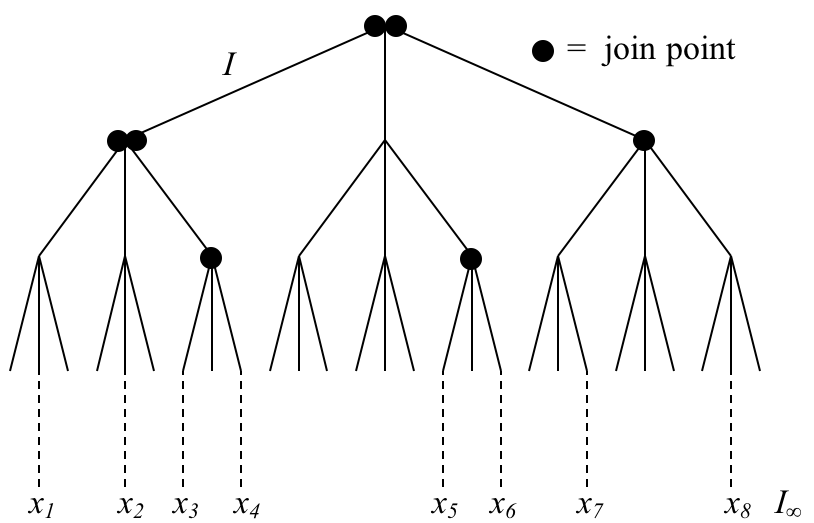}
\caption{A set of 8 points in $I_\infty$ showing the 7 join points, with  multiplicity 2  at two of the vertices}
\end{center}
\end{figure}

To establish  (\ref{intest2}) below, we will split the domain of integration into
subdomains consisting of  $n$-tuples $(\bi_1 ,\ldots,\bi_n)$ lying in different orbits
of automorphisms of the tree $I$. We will use induction over certain classes of orbit
to estimate the integrals over each such domain, with H\"{o}lder's inequality playing
a natural role at each step. It is convenient to phrase the argument using a little
terminology from group actions.

 Let $\mbox{Aut}$ be the group of automorphisms of the rooted tree $I$ that fix the
 root $\emptyset$; these automorphisms act on the infinite tree $I_\infty$ and thus
 on the $n$-tuples  $(I_\infty)^n$ in the obvious way.
  For each $n \in \bbbn$ let
 $$
S(n) = \{(\bi_1,\ldots,\bi_n) : \bi_k \in I_\infty \mbox{ for all } k =1,\ldots,n\},
 $$
for each $n \in \bbbn$ and  $l = 0,1,2,\ldots$ let
 $$
S_l(n) = \{(\bi_1,\ldots,\bi_n) : \bi_k \in I_\infty \mbox{ for all } k =1,\ldots,n |\wedge^T(\bi_1,\ldots,\bi_n)|\geq l  \},
 $$
and for each $n \in \bbbn$ and $v \in I$ let
\begin{eqnarray*}
S(\bv,n) &=& \{(\bi_1,\ldots,\bi_n) : \bi_k \in I_\infty, \bi_k \succeq \bv \mbox{ for all } k =1,\ldots,n\} \\
&=&
\{(\bi_1,\ldots,\bi_n) : \bi_k \in I_\infty, \wedge^T(\bi_1,\ldots,\bi_n)
\succeq v\mbox{ for all} k =1,\ldots,n\}.
 \end{eqnarray*}
The group Aut acts on $S(n)$ by $g(\bi_1,\ldots,\bi_n) = (g(\bi_1), \ldots, g(\bi_n))$
for $g \in \mbox{\rm Aut}$, and in particular acts on each $S_l(n)$, since if  $g \in
\mbox{Aut}$ maps $(\bi_1,\ldots,\bi_n)$ to $(\bi_1',\ldots,\bi_n')$ then $g$ maps the
join set $\wedge(\bi_1,\ldots,\bi_n)$ to the join set $\wedge(\bi_1',\ldots,\bi_n')$.
Similarly,  for each $v$, the subgroup $\mbox{\rm Aut}_v$ of $\mbox{\rm Aut}$ that fixes
the vertex $v \in I$ acts on $S(\bv,n)$. We write $\Orb_l(n)$ for the set of orbits of
$S_l(n)$ under Aut, and $\Orb(v,n)$ for the set of orbits of $S(\bv,n)$ under $\mbox{\rm Aut}_v$.
In other words,
defining an equivalence relation $\sim$ on $S(n)$ by
$$
(\bi_1,\ldots,\bi_n)  \sim (\bi_1',\ldots,\bi_n') \mbox{ if there exists } g \in
\mbox{\rm Aut}  \mbox{ such that  }   g(\bi_k) =  \bi_k'  \mbox{ for  all }  1 \leq k \leq n,
$$
then
$$\Orb_l(n)= S_l(n)/\sim  \quad  \mbox{ and } \quad  \Orb(\bv,n)= S(\bv,n)/\sim.$$

For each $l$ and $n$ we have $S_l(n) = \bigcup_{|v| = l} S(\bv,n)$ with this union disjoint.
Thus if  $|v| = l$ then the orbit $\O \in \Orb_l(n)$ restricts to an orbit $\O(v) \in \Orb(v,n)$
in the obvious way, that is
$$\O(v) =\{(\bi_1,\ldots,\bi_n) \in \O : \wedge^T(\bi_1,\ldots,\bi_n) \succeq v\}.$$

The {\it level} of a vertex $\bv \in I$ is just the length of the word $|\bv|$, and the
{\it set of levels  of a join set} is the set of $n-1$ levels of the vertices in
$\wedge(\bi_1,\ldots,\bi_n)$, counting by multiplicity. Since the set of levels of a join
set is constant across each orbit,  we may define the {\it set of levels of an orbit}
$\O \in \Orb_l(n)$ or $\O \in \Orb(\bv,n)$, written $L(\O)$,  to be the set of levels of  the join set
$\wedge(\bi_1,\ldots,\bi_n)$ of any $(\bi_1,\ldots,\bi_n) \in \O$.

Throughout this section we will be working with products over sets of levels of certain
orbits. To aid keeping track of terms, particularly when H\"{o}lder's inequality is invoked,
we use a number in square brackets above  the product sign to indicate the number of terms
in this product (for example the product in (\ref{intest}) is over $n-1$ levels).

The following proposition provides our basic estimate for the $\mu^n$-measure of $n$-tuples
of points in $I_\infty$ lying in a given orbit, in terms of the measures of cylinders at
the levels of the join classes of the orbit.

\begin{prop}\label{integerest}
Let $q>1$ and $n\geq 1$ be such that $q\geq n$.  Let $\bv \in I$ and let $\O \in \Orb(\bv,n)$.
Then
\be
\mu^n\big\{(\bi_1,\ldots,\bi_n)\in  \O \big\}
 \leq \mu(C_\bv)^{(q-n)/(q-1)} \prod_{l \in L(\O)}^{[n-1]}
 \Big(\sum_{|\bu|= l, \bu \succeq \bv} \mu(C_{\bu})^q\Big)^{1/(q-1)}. \label{intest}
\ee
\end{prop}

\noindent{\it Proof.}
We prove (\ref{intest}) by induction on $n$. When $n=1$, the only  orbit $\O\in \Orb(\bv,n)$
comprises the set of $\bi_1$ such that $\bi_1 \succeq \bv$, so
 $$\mu\big\{(\bi_1) \in \O\big\} = \mu\big\{\bi_1: \bi_1 \succeq \bv\big\} = \mu(C_\bv),$$
 which is (\ref{intest}) in this case.

Now assume inductively that  (\ref{intest}) holds for all $\O\in \Orb(\bv,n)$ for all
$\bv\in I$ and all  $1\leq n \leq n_0$. We show that  (\ref{intest}) holds with  $n= n_0+1$.

Let $ \O\in \Orb(\bv,n)$. We first consider the case where $\bv= \wedge^T(\bi_1,\ldots,\bi_n)$
for some, and therefore for all, $(\bi_1,\ldots,\bi_n)\in \O$  (thus $v$ is the top vertex
of the $\wedge (\bi_1,\ldots, \bi_n)$ in the orbit).
Then each $(\bi_1,\ldots,\bi_n)\in \O$ decomposes into $2\leq r \leq n$ subsets
$$\{(\bi_1^1,\ldots,\bi_{n_1}^1), \ldots,(\bi_1^r,\ldots,\bi_{n_r}^r) \}$$
 such that the top vertices $v_k =  \wedge^T(\bi_1^k,\ldots,\bi_{n_k}^k)$ are distinct with
 $v_k \succ v,  v_k \neq v$ and with the paths in the tree $I$ that join the $v_k$ to $v$ meeting only at $v$.
 Then $1\leq n_k <n$ for each $k$, and
 \be
n_1 +\cdots +n_r  = n. \label{weightsum}
\ee
The orbit $\O$ induces orbits $\O_k \in  \Orb(\bv,n_k)$ of  $(\bi_1^k,\ldots,\bi_{n_k}^k)$
for each $k$, and thus $\O$ may
be decomposed as a subset of the product of the $\{\O_k: 1 \leq k \leq r\}$.
Thus applying the inductive hypothesis (\ref{intest}) to each $\O_k$,
\begin{align*}
\mu^n\big\{ &(\bi_1,\ldots,\bi_n) \in \O\big\} \\
& \leq
\mu^{n_1}\big\{(\bi_1^1,\ldots,\bi_{n_1}^1) \in \O_1\big\}  \times \cdots \times
\mu^{n_r}\big\{(\bi_1^r,\ldots,\bi_{n_r}^r) \in \O_r\big\}\\
& \leq
\mu(C_{\bv})^{(q-n_1)/(q-1)} \prod_{l \in L(\O_1)}^{[n_1-1]}
 \Big(\sum_{|\bu| =l, \bu \succeq \bv} \mu(C_{\bu})^q\Big)^{1/(q-1)} \\
& \qquad\qquad \times \cdots \times
\mu(C_{\bv})^{(q-n_r)/(q-1)} \prod_{l \in L(\O_r)}^{[n_r-1]}
 \Big(\sum_{|\bu| =l, \bu \succeq \bv} \mu(C_{\bu})^q\Big)^{1/(q-1)},\\
  &= \mu(C_{\bv})^{(q-n_1- \cdots- n_r)/(q-1)} \big( \mu(C_\bv)^q\big)^{(r-1)/(q-1)}
\prod_{l \in L(\O_1)\cup \cdots \cup  L(\O_r)}^{[n_1+ \cdots+ n_r-r]}
 \Big(\sum_{|\bu| =l, \bu \succeq \bv} \mu(C_{\bu})^q\Big)^{1/(q-1)} \\
&\leq  \mu(C_{\bv})^{(q-n)/(q-1)}
\prod_{l \in L(\O)}^{[n-1]}
 \Big(\sum_{|\bu| =l, \bu \succeq \bv}  \mu(C_{\bu})^q\Big)^{1/(q-1)},
 \end{align*}
where $L(\O)$ is the complete set of levels of  ${\O}$ (including level $|v|$ with multiplicity $r-1$)  and where we have used (\ref{weightsum}).
This is (\ref{intest}) in the case where $v$ is the top vertex of the join sets of the $n$-tuples
in ${\O}$.

Finally, let $\O \in \Orb(\bv,n)$ be such that $|\wedge^T(\bi_1,\ldots,\bi_n)| = l'$ for
each $(\bi_1,\ldots,\bi_n) \in \O$ with $l' > |v|$  (so the top vertex of $\wedge(\bi_1,\ldots,\bi_n)$
is strictly below $v$). The orbit $\O$ may be decomposed into orbits $\O_\bw \in \Orb(\bw,n)$ where  $\bw \succeq \bv$ and $ |\bw| = l' $. We have shown that (\ref{intest}) holds for such $\bw$,
so summing and using H\"{o}lder's inequality,
\begin{align*}
\mu^n\big\{(\bi_1,& \ldots,\bi_{n}) \in \O\big\}
\\
& \leq\sum_{|\bw| = l', \bw \succeq \bv} \bigg(  \mu(C_\bw)^{(q-n)/(q-1)} \prod_{l \in L(\O_\bw)}^{[n-1]}
 \Big(\sum_{|\bu| =l, \bu \succeq \bw} \mu(C_{\bu})^q\Big)^{1/(q-1)}\bigg)\\
& \leq \Big(\sum_{|\bw| = l', \bw \succeq \bv}   \mu(C_\bw)\Big)^{(q-n)/(q-1)} \prod_{l \in L(\O)}^{[n-1]}
 \Big( \sum_{|\bw| = l', \bw \succeq \bv} \Big( \sum_{|\bu| = l, \bu \succeq \bw}  \mu(C_{\bu})^q \Big)\Big)^{1/(q-1)}\nonumber\\
 & = \mu(C_\bv)^{(q-n)/(q-1)} \prod_{l \in L(\O)}^{[n-1]}
 \Big(  \sum_{|\bu| = l, \bu \succeq \bv}  \mu(C_{\bu})^q\Big)^{1/(q-1)}.
  \end{align*}
This completes the induction and the proof.
$\Box$

\medskip

Proposition \ref{integerest} would be enough for our purposes when $q$ is an integer.
However, when estimating (\ref{intest2}) for  a non-integer  $q>1$ we need a generalization
where one of the points $\bj \in I_{\infty}$ is distinguished.
The proof of Proposition \ref{fracest}  again uses induction on join sets and H\"{o}lder's
inequality, but the argument is more intricate than that of Proposition \ref{integerest}
on which it depends.

Let $ 1\leq p \leq n$,  $0 \leq l_1 < l_2 < \ldots <l_p$ be  levels and $m_1, \ldots, m_p
\in \bbbn $ be such that $m _1+\cdots + m_p = n$. For each  $\bj \in I_{\infty}$ and
$1 \leq r \leq p$, write  $\bj_r =\bj|_{l_r}$ and, given  an orbit $\O_{r} \in \Orb_{l_r}(m_r)$,
write $\O_{r}(\bj_r)  \in \Orb(\bj_r,m_r)$
for the suborbit under automorphisms fixing $\bj_r$. 

In the next proposition we integrate over $y$ powers of the measures of those $(\bi_1 \ldots,\bi_{n})$ for which the joins of the distinguished point $y$ with the $\bi_k$ are  the vertices  $\bj_r=\bj|_{l_r}$ lying on the path from $\emptyset$ to $y$ in the tree $I_\infty$ at levels $l_r$, and which, for each $r$, the set of $\bi_k$ such that $\bi_k \wedge \bj = \bj_r$ lie in a given  orbit fixing  $\bj_r$. Thus the set of join levels of $(\bi_1 \ldots,\bi_{n},\bj)$ comprises the levels $l_r$ of joins with $\bj$ together with the join levels of each of the orbits 
$\O_{r}(y_r)$, which are just the join levels  $L(\O_r)$ of the orbits $\O_{r}$. Thus the integral in the next proposition  is bounded by a product of sums over taken over join levels of these two types: $ L=\{l_1,\ldots,l_p,L(\O_1),\ldots, L(\O_p)\}$.

\begin{prop}\label{fracest}
Let $q>1$ and let $n$ be an integer with $n\geq q-1$.
As above, let $ 1\leq p \leq n$, let  $0 \leq l_1 < l_2 < \ldots <l_p$ be  levels and
let $m _1+\cdots + m_p = n$. For each $r=1,\ldots,p$ let $\O_r \in \Orb_{l_r}(m_r)$ be given.
Then
\begin{align}
\int_{\bj \in I_\infty}
\mu^n\big\{(\bi_1 \ldots,\bi_{n}):
(\bi_{1},\ldots,\bi_{m_{1}})&  \in \O_{1}(\bj_1), \ldots ,
 (\bi_{n-m_p +1},\ldots,\bi_{n}) \in \O_{p}(\bj_p)\big\}
^{(q-1)/n} d\mu(\bj) \nonumber\\
&\leq
\prod_{l \in L}^{[n]}
 \Big(\sum_{|\bu| =l} \mu(C_{\bu})^q\Big)^{1/ n}, \label{qintest}
\end{align}
where $ L$ denotes the aggregate set of levels of $\{l_1,\ldots,l_p,L(\O_1),\ldots, L(\O_p)\}$.
\end{prop}

\noindent{\it Proof.}
We proceed by induction on $r \leq p$, starting with $r= p$ and working backwards to $r=1$,
taking as the inductive hypothesis:

\noindent For all $\bj_r \in I_{ l_r}$,
\begin{align}
\int_{\bj \succeq \bj_r}
\mu^{m_r+ \cdots + m_p}
\big\{(&\bi^{r}_{1},\ldots,\bi^{r}_{m_{r}},\ldots,
\bi^{p}_{1},\ldots,\bi^{p}_{m_{p}})\nonumber \\
&: (\bi^{r}_{1},\ldots,\bi^{r}_{m_{r}}) \in \O_{r}(\bj_{r}),\ldots,
(\bi^{p}_{1},\ldots,\bi^{p}_{m_{p}}) \in \O_{p}(\bj_p)\big\}
^{(q-1)/n} d\mu(\bj) \nonumber\\
\leq \mu(C_{\bj_r}&)^{(n - n_r )/n}
\prod_{l \in L_r}^{[ n_r]}
 \Big(\sum_{|\bu| =l, \bu \succeq \bj_r }  \mu(C_{\bu})^q\Big)^{1/n},\label{indhop}
\end{align}
where $n_r = m_r+\cdots+m_p$ and  $ L_r$ denotes the set of levels of $\{l_r,\ldots,l_p,L(\O_r),
\ldots, L(\O_p)\}$ counted by multiplicity (so that  $ L_r$ consists of $m_r+\cdots+m_p = n_r$ levels).

To start the induction, we apply Proposition \ref{integerest} to get, for each $\bj_p \in I_{ l_p}$,

\begin{align}
\int_{\bj \succeq \bj_p}\mu^{m_p}
\big\{&(\bi^{p}_{1},\ldots,\bi^{p}_{m_{p}}) \in \O_{p}(\bj_p)\big\}
^{(q-1)/n}  d\mu(\bj)\label{exp1} \\
&
\leq  \int_{\bj \succeq \bj_p}\bigg[
\mu(C_{\bj_p}) ^{(q - m_p)/(q-1)}
\prod_{l \in L(\O_p(\bj_p))}^{[ m_p-1]}
 \Big(\sum_{|\bu| =l, \bu \succeq \bj_p}\mu(C_{\bu})^q
\Big)^{1/(q-1)}\bigg]^{(q-1)/n} d\mu(\bj)  \nonumber\\
&
=
\mu(C_{\bj_p})^{(n - m_p)/n} \big(\mu(C_{\bj_p})^{q}\big)^{1/n}
\prod_{l \in  L(\O_p)}^{[ m_p-1]}
 \Big(\sum_{|\bu| =l, \bu \succeq \bj_p}  \mu(C_{\bu})^q
\Big)^{1/n}  \nonumber \\
& =
\mu(C_{\bj_p}) ^{(n - n_p)/n}
\prod_{l \in L_r}^{[ n_p]}
 \Big(\sum_{|\bu| =l, \bu \succeq \bj_p } \mu(C_{\bu})^q\Big)^{1/n},  \nonumber
 \end{align}
on incorporating $(\mu(C_{\bj_p})^{q})^{1/n}$ in the main product and noting that $m_p = n_p$.
(Observe that this remains valid if $m_p =1$, in which case  the measure in (\ref{exp1}) is
at most $\mu(C_{\bj_p})$.) This establishes the inductive hypothesis (\ref{indhop}) when $r=p$.

Now assume that (\ref{indhop}) is valid for $r= k,\ldots, p$ for some $ 2\leq k \leq p$. Then
for each $\bj_{k-1} \in I_{l_{k-1}}$,
\begin{align}
I &:= \int_{\bj \succeq \bj_{k-1}}
\mu^{m_{k-1} + m_k + \cdots + m_p}
\big\{(\bi^{k-1}_{1},\ldots,\bi^{k-1}_{m_{k-1}},\bi^{k}_{1},\ldots,\bi^{k}_{m_{k}} ,\ldots,
\bi^{p}_{1},\ldots,\bi^{p}_{m_{p}})\nonumber \\
&\qquad\qquad\qquad \qquad\qquad\qquad: (\bi^{k-1}_{1},\ldots,\bi^{k-1}_{m_{k-1}}) \in \O_{k-1}(\bj_{k-1}),
(\bi^{k}_{1},\ldots,\bi^{k}_{m_{k}}) \in \O_{k}(\bj_{k}),\nonumber \\
&\qquad\qquad\qquad \qquad\qquad\qquad\qquad\qquad\qquad\qquad \ldots,
(\bi^{p}_{1},\ldots,\bi^{p}_{m_{p}}) \in \O_{p}(\bj_p)\big\}
^{(q-1)/n}
d\mu(\bj)  \nonumber \\
&\leq \mu^{m_{k-1}} \big\{
(\bi^{k-1}_{1},\ldots,\bi^{k-1}_{m_{k-1}}) \in \O_{k-1}(\bj_{k-1})\big\}
^{(q-1)/n}\nonumber \\
&\qquad\qquad\times  \sum_ {\bj_k \succeq \bj_{k-1}, |\bj_k| = l_k}
\int_{\bj \succeq \bj_{k}}
\mu^{m_k+ \cdots + m_p}
\big\{(\bi^{k}_{1},\ldots,\bi^{k}_{m_{k}},\ldots,
\bi^{p}_{1},\ldots,\bi^{p}_{m_{p}})\nonumber \\
&\qquad\qquad\qquad \qquad: (\bi^{k}_{1},\ldots,\bi^{k}_{m_{k }}) \in \O_{k}(\bj_{k}),\ldots,
(\bi^{p}_{1},\ldots,\bi^{p}_{m_{p}}) \in \O_{p}(\bj_p)\big\}
^{(q-1)/n}
 d\mu(\bj)  \nonumber\\
  & \leq
 \bigg[\mu(C_{\bj_{k-1}})^{(q - m_{k-1})/n}
\prod_{l \in L(\O_{k-1} (\bj_{k-1}))}^{[ m_{k-1}-1]}
 \Big(\sum_{|\bu| =l, \bu \succeq \bj_{k-1}}  \mu(C_{\bu})^q
\Big)^{1/n}\bigg]\nonumber\\
&\qquad \qquad \qquad \qquad\qquad\times \sum_ {{\bj_k \succeq \bj_{k-1}, |\bj_k| = l_k}}
\bigg(\mu(C_{\bj_k}) ^{(n - n_k)/n}
\prod_{l \in L_k}^{[ n_k]}
 \Big(\sum_{|\bu| =l, \bu \succeq \bj_k }  \mu(C_{\bu})^q\Big)^{1/n}\bigg),\label{jaykay}
 \end{align}
 where in obtaining the last inequality we have used Proposition \ref{integerest} to
 estimate the first part and the inductive hypothesis (\ref{indhop}) for the second part.
Using H\"{o}lder's inequality for each $\bj_{k-1}$:
\begin{align*}
\sum_ {{\bj_k \succeq \bj_{k-1}, |\bj_k| = l_k}}&\bigg(\mu(C_{\bj_k})^{(n - n_k )/n}
\prod_{l \in L_k}^{[ n_k ]}
 \Big(\sum_{|\bu| =l, \bu \succeq \bj_k }  \mu(C_{\bu})^q\Big)^{1/n} \bigg)\\
& \leq\Big(\sum_ {{\bj_k \succeq \bj_{k-1}, |\bj_k| = l_k}}\mu(C_{\bj_k})\Big)^{(n - n_k )/n}
\prod_{l \in L_k}^{[ n_k ]}
 \Big(\sum_ {{\bj_k \succeq \bj_{k-1}, |\bj_k| = l_k}} \Big(\sum_{|\bu| =l, \bu \succeq \bj_k } \mu(C_{\bu})^q \Big)\Big)^{1/n} \\
& = \mu(C_{\bj_{k-1}})^{(n - n_k )/n}
\prod_{l \in L_k}^{[ n_k ]}
 \Big(\sum_{|\bu| =l, \bu \succeq \bj_{k-1} } \mu(C_{\bu})^q\Big)^{1/n}.
\end{align*}
Thus from (\ref{jaykay})
 \begin{align*}
I \leq \mu &(C_{\bj_{k-1}})^{(n - n_k - m_{k-1} )/n} \big( \mu(C_{\bj_{k-1}})^q\big)^{1/n}
 \\
&\times \prod_{l \in L(\O_{k-1} (\bj_{k-1}))}^{[ m_{k-1}-1]}
 \Big(\sum_{|\bu| =l, \bu \succeq \bj_{k-1}}  \mu(C_{\bu})^q
\Big)^{1/n}
\prod_{l \in L_k}^{[ n_k]}
 \Big(\sum_{|\bu| =l, \bu \succeq \bj_{k-1} } \mu(C_{\bu})^q\Big)^{1/n}\\
 = \mu &(C_{\bj_{k-1}})^{(n - n_k - m_{k-1} )/n}
\prod_{l \in L_{k-1}}^{[ m_{k-1}+n_k]}
 \Big(\sum_{|\bu| =l, \bu \succeq \bj_{k-1} } \mu(C_{\bu})^q\Big)^{1/n} ,
 \end{align*}
which is (\ref{indhop}) with $r=k-1$, noting that $m_{k-1}+n_k = n_{k-1}$.

Finally, taking $r=1$ in (\ref{indhop}) and noting that $n_1 = n$,
\begin{align*}
\int_{\bj \succeq \bj_1 }
\mu^{n}
\big\{(\bi_{1},\ldots,\bi_{n}): (\bi_{1},\ldots,\bi_{m_{1}}) \in \O_{1}(\bj_{1}),&\ldots,
(\bi_{n-n_p+1},\ldots,\bi_{n}) \in \O_{p}(\bj_p)\big\}
^{(q-1)/n}
d\mu(\bj)  \\
&  \leq
\prod_{l \in L_1}^{[ n]}
 \Big(\sum_{|\bu| =l, \bu \succeq \bj_1} \mu(C_{\bu})^q\Big)^{1/n},
\end{align*}
and summing  over all $\bj_1$ at level $l_1$ and using H\"{o}lder's inequality again, gives (\ref
{qintest}).
$\Box$

\medskip

To use Proposition \ref{fracest} to determine when the integral in (\ref{intest2}) converges
we need to bound the number of orbits that have prescribed sets of join levels. Let
$0 \leq k_1\leq \cdots \leq k_n$ be (not necessarily distinct) levels. Two types of join levels arise in Proposition  \ref{fracest}: those on a distinguished path at levels  $l_r$ and the  join levels  $\O_r$ of the subsidiary orbits, see the remark before Proposition \ref{fracest}.
\begin{align}
N(k_1,\ldots, k_n) = \#\Big\{(l_1,&\ldots,l_p, \O_1,\ldots,\O_p): 1 \leq p \leq n,
0 \leq l_1<\cdots <l_p, \nonumber\\
&\O_r \in \Orb_{l_r} (m_r)\mbox{ for some } m_r \mbox{ where }
m_1 + \cdots + m_r = n, \nonumber\\
&\mbox{ such that }
\{l_1,\ldots,l_p,L(\O_1),\ldots,L(\O_p)\} = \{k_1,\ldots, k_n\} \Big\}. \label{numbdef}
\end{align}

\begin{lem}\label{count}
Let  $n \in \bbbn$ and $0<\lambda<1$. Then
$$\sum_{0 \leq k_1\leq \cdots \leq k_n}N(k_1,\ldots,k_n)  \lambda^{(k_1+\cdots +k_n)/n} <\infty.
$$
\end{lem}
{\it Proof.}
The crucial observation here is that we may find an upper bound for $N(k_1,\ldots,k_n)$ that depends on $n$ but not on the particular levels $(k_1,\ldots, k_n)$.

Let  $N_0 (k_1,\ldots, k_n)$ be the total number of orbits in $\Orb(\emptyset,n+1)$
(where $\emptyset$ is the root of the tree $I$) with levels $0 \leq k_1\leq \cdots \leq k_n$.
 Every join set  with levels   $0 \leq k_1\leq \cdots \leq k_n$ may be obtained by adding
a vertex  at level $k_{n} $ of the form $ x_i|_{k_{n}}$ to a join set $\wedge(x_1,\ldots,x_{n})$ with levels
$0 \leq k_1\leq \cdots \leq k_{n-1}$ for some $1 \leq i \leq n$, and this may be done in at most $n$
ways. It follows that
$N_0 (k_1,\ldots, k_{n}) \leq n  N_0 (k_1,\ldots, k_{n-1})$, so since $N_0 (k_1)=1$, we obtain
$N_0 (k_1,\ldots, k_n)  \leq n!$.

The number $N(k_1,\ldots, k_n)$ given by (\ref{numbdef}) is no more than the number of orbits in $\Orb(\emptyset,n+1)$ having levels $0 \leq k_1\leq \cdots \leq k_n$ with a subset of the join vertices (to within equivalence) of each member of the orbit distinguished to correspond to levels $l_1,\ldots,l_p$.
But given $(x_1,\ldots,x_{n+1}) \in \O$ where $\O$ has join levels $k_1,\ldots, k_n$, there are at most $2^n$ ways of choosing a distinguished subset $(y_1,\ldots,y_p)$ of the  join set $\wedge(x_1,\ldots,x_{n+1})$;  these vertices then determine $m_r$ as well as
$\O_r(y_r) = \O_r \in \Orb_{l_r} (m_r)$ for $r=1,\ldots,p$.
(Note that we are considerably over-counting since the same contributions to (\ref{numbdef}) may come from different orbits of  $\Orb(\emptyset,n+1)$.)   Hence
$$N(k_1,\ldots, k_n) \leq 2^n N_0 (k_1,\ldots, k_n) \leq 2^n n! .$$
 Thus
 \begin{align}
\sum_{0 \leq k_1\leq \cdots \leq k_n}N(k_1,\ldots,k_n)  \lambda^{(k_1+\cdots +k_n)/n}
&\leq
2^n n!\sum_{0 \leq k_1\leq \cdots \leq k_n}  \lambda^{( k_1+\cdots +k_n)/n}\nonumber\\
&\leq
2^n n! \sum_{k=0}^{\infty}P(k)  \lambda^{k/n},\label{series}
\end{align}
where $P(k)$ is the number of distinct ways of partitioning the integer $k$ into a sum
of $n$ integers $k = k_1 + \cdots + k_{n}$ with $0 \leq k_1\leq \cdots \leq k_n$.
Since $P(k)$ is polynomially bounded (trivially $P(k) \leq (k+1)^{n-1}$), (\ref{series})
converges for $0< \lambda<1$.
$\Box$

\medskip
To obtain the main estimate, we use Lemma \ref{count} to count the domains of integration to which we apply Proposition \ref{fracest}.

Let $f: \bbbn_0 \to \bbbr_+$ be a function.
Define the {\it multipotential kernel} $\phi: I^{n+1} \to \bbbr_+$ to be the product of
$f$ evaluated at the levels of the vertices of each join set, that is
\be
 \phi(\bi_1,\ldots,\bi_n,\bj) = f(l_1)f(l_2)\cdots f(l_n) \mbox{ where }
 L(\wedge(\bi_1,\ldots,\bi_n,\bj)) = \{l_1,l_2,\ldots l_n\}. \label{defphi}
\ee

\begin{theo}\label{mainint}
Let $n\in \bbbn$ and $q>1$ with $n\leq q< n+1$. Suppose that
\begin{equation}
\limsup_{l \to \infty}  \frac{\log\big( f(l)^{q-1}\sum_{|\bu| =l} \mu(C_{\bu})^q\big)}
{  l} < 0.\label{condition}
\end{equation}
Then, with $\phi$ as in $(\ref{defphi})$,
\begin{equation}
J := \int\bigg[
 \int \cdots \int
 \phi(\bi_1,\ldots,\bi_n,\bj) d\mu(\bi_1)\ldots d\mu(\bi_n)\bigg]^{(q-1)/n} d\mu(\bj)<\infty.
 \label{intest2}
 \end{equation}
\end{theo}
{\it Proof.}
For each $\bj \in I_\infty$ we decompose the integral inside the square brackets as a sum
of integrals taken over all $p$, all $0 \leq l_1<\cdots<l_p$, all  $m_1,\dots,m_p\geq 1$
such that $m_1+\cdots+m_p = n$, and all orbits
$\O_1\in \Orb_{l_1}(m_1),\ldots,\O_p\in  \Orb_{l_p}(m_p)$. As before, for each $r$ we write
$\bj_r = \bj|_{l_r}$,  and  $\O_r(y_r)$ for the suborbit of $\O_r$ in $\Orb(y_r,m_p)$. Thus,
using the power-sum inequality, (\ref{defphi}) and (\ref{qintest}), and noting that $\phi$ depends only on the levels of the join sets,

\begin{align*}
&J = \\
&\int \bigg[\hspace{-0.3cm} \sum_{{\scriptsize
\begin{array}{c}
0 \leq l_1<\cdots<l_p \\
m_1+\cdots+m_p = n   \\
\O_1,\ldots,\O_p
\end{array}}}
\hspace{-0.3cm} \int_{(\bi_{1},\ldots,\bi_{m_{1}}) \in \O_{1}(y_1)}\hspace{-0.5cm} \cdots
\int_{(\bi_{n-m_p +1},\ldots,\bi_{n}) \in \O_{p}(y_p)}
\hspace{-1.5cm} \phi(\bi_1,\ldots,\bi_n,\bj) d\mu(\bi_1)\ldots d\mu(\bi_n)\bigg]^{(q-1)/n}
\hspace{-0.5cm} d\mu(\bj)  \\
&
\leq\hspace{-0.5cm}\sum_{{\scriptsize
\begin{array}{c}
0 \leq l_1<\cdots<l_p \\
m_1+\cdots+m_p = n   \\
\O_1,\ldots,\O_p
\end{array}}}
\int \bigg[\mu^n \big\{(\bi_{1},\ldots,\bi_{n}): (\bi_{1},\ldots,\bi_{m_{1}}) \in \O_{1}(y_1),\\
&\hspace{4cm}\ldots,
(\bi_{n-m_p +1},\ldots,\bi_{n}) \in \O_{p}(y_p)\big\}
\phi(\bi_1,\ldots,\bi_n,\bj)\bigg]^{(q-1)/n} d\mu(\bj)  \\
&
\leq\sum_{{\scriptsize
\begin{array}{c}
0 \leq l_1<\cdots<l_p \\
m_1+\cdots+m_p = n   \\
\O_1,\ldots,\O_p
	\end{array}}}
 \prod_{l \in L}^{[n]}
 \Big(f(l)^{q-1} \sum_{|\bu| = l}  \mu(C_{\bu})^q\Big)^{1/ n},
\end{align*}
where the products are over the set of levels
$L = \{l_1,\ldots,l_p,L(\O_1),\ldots, L(\O_p)\}$ counted with repetitions
Condition (\ref{condition}) implies that
$f(l)^{q-1}\sum_{|\bu| =l} \mu(C_{\bu})^q \leq c \lambda^l$ for all $l$, for some
$c>0$ and some $ \lambda <1$. Thus, with $N(k_1,\ldots,k_n)$ as in (\ref{numbdef}),
\begin{align*}
I &
\leq\sum_{0 \leq k_1 \leq \ldots \leq k_n}N(k_1,\ldots,k_n)
 \prod_{i=1}^{n}
 \Big(f(k_i)^{q-1} \sum_{|\bu| = k_i}  \mu(C_{\bu})^q\Big)^{1/ n}\\
&\leq\sum_{0 \leq k_1 \leq \ldots \leq k_n}
 N(k_1,\ldots,k_n)\prod_{i=1}^{n}
 \Big(c\lambda^{k_i})^{1/ n}\\
&
\leq
 c\sum_{0 \leq k_1 \leq \ldots \leq k_n}N(k_1,\ldots,k_n)  \lambda^{(k_1+\cdots +k_n)/n} <\infty
 \end{align*}
using Lemma \ref{count}.
$\Box$

\section{Proofs of main results}\label{ranpro}
\setcounter{equation}{0}
\setcounter{theo}{0}

We can now complete the proof of our main results.

Firstly, to enable us to work separately with upper and lower indices, it is convenient to have a variant  of (\ref{Holder}) under a weakened  H\"{o}lder condition

\begin{lem}\label{Lem:Holder}
Suppose that $f:\bbbr^N \to \bbbr^d$ and that for some compact interval $K \subseteq \bbbr^N$ and  $0 < \alpha \le 1$,
there exist a sequence  $r_n \searrow 0$ and a constant $c$ such that
\begin{equation}\label{Eq:p-holder}
\sup_{x, y \in K: |x-y|\le r_n} \big|f(x) - f(y) \big|\le c\, r_n^\alpha.
\end{equation}
Then for all $q>0, q \neq 1$ and every finite Borel measure $\mu$ with support contained in $K$,
we have
\[
\lqd ( \mu_f)   \leq \min\Big\{d,\, \frac 1 \alpha\, \uqd (\mu) \Big\}.
\]
In particular, if the
generalized $q$-dimension $D_q (\mu)$ exists, then
\[
\lqd ( \mu_f)   \leq \min\Big\{d,\, \frac 1 \alpha\, D_q (\mu) \Big\}.
\]
\end{lem}
\noindent{\it Proof.}\, Let $\rho_n = c\, r_n^\alpha$ for all $n \ge 1$.
Then (\ref{Eq:p-holder}) implies that
\[
\begin{split}
\int_{\bbbr^d} \mu_f\big(B(z, \rho_n)\big)^{q-1} d\mu_f(z) &=
\int_{K} \mu_f\big(B(f(x), \rho_n)\big)^{q-1} d\mu(x)\\
& \ge \int_{K} \mu\big(B(x, r_n)\big)^{q-1} d\mu(x)\\
\end{split}
\]
for $q>1$, with the reverse inequality for $0<q<1$.
Hence, in both cases,
\[
\begin{split}
\liminf_{n \to \infty} \frac{\int_{\bbbr^d} \mu_f\big(B(z, \rho_n)\big)^{q-1}
d\mu_f(z)}{(q-1)\log \rho_n} &\le \limsup_{n \to \infty}
\frac{\int \mu\big(B(x, r_n)\big)^{q-1} d\mu(x)} {(q-1)\log \rho_n}\\
&= \frac 1 \alpha\, \uqd (\mu).
\end{split}
\]
$\Box$

The following well-known lemma on the modulus of continuity of a
Gaussian process follows from  \cite[Corollary 2.3]{Dudley73}.

\begin{lem}\label{Lem:Modulus}
Let $X_0: \R^N \to \bbbr$ be a centered Gaussian random field which satisfies Condition
{\bf (C1)} for some $\psi$ with $0 < \a_* \le \a^* < 1$. Given a compact interval
$K \subseteq \R^N$, let
$$ \omega_{X_0}(\de) = \sup_{\mbox{\footnotesize {$\begin{array} {c} x, \ x + y \in K\\
|y| \le \de \end{array}$}}} |X_0(x+y) - X_0(x)| $$ be the uniform
modulus of continuity of  $X_0(x)$ on $K$. Then there exists a finite constant $C_{4} > 0$
such that
$$
\limsup_{\de \to 0} \frac{\omega_{X_0} (\de)} { \sqrt{ \psi(\de) \log
\frac1 {\de} }}  \le C_{4}, \quad \hbox{a.s.}
$$
\end{lem}

Assuming Condition {\bf (C1)} we can now obtain  a.s. upper
bounds for the generalized $q$-dimensions of the image measure $\mu_X$ of $\mu$.
\begin{prop}\label{Prop:UP}
Let $X: \R^N\to \bbbr^d$ be the Gaussian random field defined by (\ref{def:X}) such that the
associated random field $X_0$ satisfies Condition {\bf (C1)}, and let $q > 0,q\neq 1$.
Let $\mu$ be a Borel probability measure $\mu$ on $\R^N$ with compact support.
\begin{itemize}
\item[(i)]\, If $ 0<\a_* = \a^*=\alpha < 1$, then
$$
\lqd ( \mu_X)   \leq \min\Big\{d,\, \frac 1 \alpha\, \lqd (\mu) \Big\}  \,\, \hbox{a.s.}\quad
\mbox{ and } \quad
\uqd ( \mu_X)   \leq \min\Big\{d,\, \frac 1 \alpha\, \uqd (\mu) \Big\} \,\, \hbox{a.s.}
$$

\item[(ii)]\, If\, $0 < \a_* \leq \a^* < 1$, then
$$
\lqd (\mu_X)   \leq \min\Big\{d,\, \frac 1 {\alpha^*}\, \uqd (\mu) \Big\}  \,\, \hbox{a.s.}
\quad \mbox{ and } \quad
\uqd (\mu_X)   \leq \min\Big\{d,\, \frac 1 {\alpha_*}\, \uqd (\mu) \Big\} \,\, \hbox{a.s.}
$$
\end{itemize}
\end{prop}
{\it Proof.}
This follows from combining (\ref{Holder}) and Lemmas \ref{Lem:Holder} and \ref{Lem:Modulus}. If
$0<\beta'<\beta <\alpha^*$ then from (\ref{Def:phiup}) $\psi(r) \leq c r^{2\beta}$ for all small $r$. By Lemma \ref{Lem:Modulus}  there is a.s. a random $C$ such that $\omega_{X_0}(\delta) \leq C\delta^{\beta'}$ for all small $\delta$, so that 
$|X_0(x) - X_0(y)|  \leq C\delta^{\beta'}$ for all  $x$ and $y$  in the support of $\mu$ with  $|x-y|$ sufficiently small. By Lemma \ref{Lem:Holder} $\lqd (\mu_X)   \leq  \frac 1 {\beta'}\, \uqd (\mu)$; this is true for all 
$0<\beta' <\alpha^*$, giving the first inequality of (ii). The other inequalities are derived in a similar way.
$\Box$

As is often the case when finding dimensions or generalized dimensions,
lower bounds are more elusive than upper bounds.
For this, we use the strong local $\psi$-nondeterminism of $X$
along with Corollary \ref{corums}  to bound the probability that points $X(x_1),
\ldots,X(x_n)$ all lie in  the ball $B(X(y),r)$ in terms of a multipotential
kernel defined using an ultrametric $d_a$.  We then apply
Theorem \ref{mainint}  to bound an integral involving this kernel.

Recall that an isotropic multivariate Gaussian random variable $Z$ in $\bbbr^d$
with variance $\sigma^2$ satisfies
$$
\P\big\{ |Z-u| \leq r \big\}
\leq  c \Big( \frac{r}{\sigma} \Big)^s,  \qquad u \in  \bbbr^d, r>0
$$
for all $0<s \leq d$, for some constant $c\equiv c_{d,s}$. Since the conditional
distributions in a Gaussian process are still Gaussian, it follows from
the strong local nondeterminism of $X_0$  given by {\bf (C2)} that, for $t_0>0$, there exists $c>0$
such that
\begin{equation}
\P\big\{ |X(x)-u| \leq r \;\big | \;X(y) \, :\, t \leq |x-y| \leq t_0\big\} \leq
cr^s  \psi(t)^{- s/2} \label{lnd2x}
\end{equation}
for all $u \in  \bbbr^d$, $x \in \bbbr^N$, $t>0$ and $ r>0$.

We now introduce the multipotential kernel that will be used when applying Theorem
\ref{mainint}. Let $m \geq 2$ and  $a  \in {\cal A}_m$, see (\ref{translates}).  Recall that ${\cal C}_{k}^a$
denotes the set of $k$th level cubes in the hierarchy of $m$-ary half-open cubes that define
the metric $d_a$, i.e. the cubes of ${\cal C}_{k}$ translated by the vector $-a$,
see (\ref{cubes}) and (\ref{cka}). Let   ${\cal C}^a = \cup_{k=0}^\infty {\cal C}_{k}^a$.

For a given $a  \in {\cal A}_m$ and  $w,z \in [0,\frac{1}{2})^N$, let $w\wedge z$
denote the smallest cube $C \in {\cal C}^a$ such that $w,z \in C$.
Then if $w_1,\ldots, w_{n+1}$ are distinct points of $[0,\frac{1}{2})^N$ there is a
uniquely defined set of $n$ {\em join cubes}, $C_1,\ldots,C_n \in {\cal C}^a$, with
the property that $w_i\wedge w_j  $ is  one of the cubes $C_l$ for all $i\neq j$.
To ensure that a set of $n+1$ points has exactly $n$ join cubes we regard a join cube
$C$ as having {\it multiplicity} $r \geq 1$ if $r$ is the greatest integer such that
there are distinct $w_{j_1},\ldots, w_{j_{r+1}}$ with $w_{j_p}\wedge w_{j_q}$ as the
cube $C$ for all $1 \leq p<q \leq r+1$. 

The hierarchy of cubes ${\cal C}^a$ may be identified in a natural way with the vertices of the $m^N$-ary tree $I$ of  Section \ref{integralest}, with  the points of
$[0,1)^N - a$ identified with the infinite tree $I_\infty$, and with corresponding join points identified. This identification provides an isometry between $[0,1)^N - a$ under the ultrametric of Section 3 and the corresponding tree of Section \ref{integralest} under the natural metric.

For a simple example, taking $m=2$ and $N=1$,  ${\cal C}_{0}^a$ is the interval $[-a,-a+1]$ and ${\cal C}_{k}^a$ consists of $2^k$ half-open intervals of length  $2^{-k}$ each of which contains $2$ disjoint intervals  of  ${\cal C}_{k+1}^a$. 
This hierarchy of intervals is indexed by the binary tree $I= \cup_{k=0}^\infty\{1,2\}^k$ under the identification 
$$\textstyle h(i_1,\ldots,i_k) = [ \sum_{l=1}^k (i_l -1)2^{-l} -a , \  \sum_{l=1}^k (i_l -1)2^{-l}+ 2^{-k}-a)$$ 
with the infinite words of $I_\infty = \{1,2\}^{\mathbb{N}}$ identified with $[-a,-a+1)$ by $h(i_1,i_2,\ldots) =  \sum_{l=1}^\infty (i_l -1)2^{-l}-a$. In particular, for $x_1, x_2 \in I_\infty$, the interval  $h(x_1 \wedge x_2)=h(x_1) \wedge h(x_2)$ is the smallest interval in ${\cal C}^a$ containing both of the points $h(x_1)$ and $h (x_2)$. Henceforth we will identify the tree with the  hierarchy of cubes in this way, referring to whichever representation is most convenient.

Write $k(C)$ for the  {\it level} of the cube $C\in{\cal C}^a$, so that $C \in {\cal C}_{k(C)}^a$.
We define the {\em multipotential kernel} $\phi_a$ by
\be
\phi_a(w_1,\ldots, w_{n+1}) = m^{k(C_1)}m^{k(C_2)}\cdots m^{k(C_n)}, \label{multpotdef}
\ee
where $C_1,\ldots,C_n $ are the join cubes of $w_1,\ldots, w_{n+1}$ in the hierarchy ${\cal C}^a$.

The following proposition uses strong local nondeterminism, in the form of (\ref{lnd2x}),
to estimate inductively the probability that the images of a set of points all lie inside a ball.

\begin{prop}\label{lnd}
Let $X: \R^N\to \R^d$ be the Gaussian random field defined by $(\ref{def:X})$, and assume the
associated random field $X_0: \R^N\to \R$ satisfies Condition {\bf (C2)}.
Given $N,d,n,s$ and $ \alpha$, where   $0<s \leq d$ and $\alpha > \alpha^*$, there are
positive constants $c_2$ and $r_0$ and an integer $m \geq 2$ such that, for all $x_1,x_2,
\ldots, x_n,y \in [0,\frac{1}{2})^N$, we may choose a vector $a \in {\cal A}_m$, see
{\rm (\ref{translates})}, such that for all $0<r \leq r_0$
\begin{align}
\P\Big\{
 |X(y)-X(x_1)|  \leq r, |X(y)-X(x_2)|  \leq r,&\ldots,|X(y)-X(x_n)| \leq r\Big\}\nonumber\\
& \leq  c_2 r^{sn}\phi_a(x_1,x_2,\ldots,x_n,y)^{\alpha s}.
\label{probdists4}
\end{align}
In particular, for all $x_1,x_2,\ldots, x_n,y \in [0,\frac{1}{2})^N$ and $0<r \leq r_0$,
\begin{align}
\P\Big\{
 |X(y)-X(x_1)|  \leq r, |X(y)-X(x_2)|  \leq r, & \ldots,|X(y)-X(x_n)| \leq r\Big\}\nonumber\\
& \leq  c_2 r^{sn}\sum_{a{ \in \cal A}_m} \phi_a(x_1,x_2,\ldots,x_n,y)^{\alpha s}.
\label{probdists5}
\end{align}
\end{prop}

\noindent{\it Proof.}
Let $m= 2n^2N+2$ and $c_0 = \max \{8m(m-1), N^{1/2}\}$. By Corollary \ref{corums},
given   $x_1,x_2,\ldots, x_n,y \in [0,\frac{1}{2})^N$, there exists  $a \in {\cal A}_m$
such that
\be
c_0^{-1}|z-w| \leq  d_a(z,w)  \leq c_0 |z-w|, \qquad z,w \in \{ x_1,x_2,\ldots, x_n,y\};
\label{equiv1}
\ee
thus $d_a$ restricted to the set of points $ \{ x_1,x_2,\ldots, x_n,y\}$ is
equivalent to the Euclidean metric with constant $c_0$.

We now appeal to the strong local $\psi$-nondeterminism of $X$.
For $i=2,3,\ldots,n$ let $w_i$ be the point (or one of the points) from $\{x_1,
\ldots,x_{i-1},y\}$ such that $d_a(w_i,x_i)$ is least. By local nondeterminism (\ref{lnd2x}),
noting the equivalence of the metrics (\ref{equiv1}), there are constants $c_1$ and $r_0$
such that
\begin{equation}\label{Eq:conp}
\begin{split}
\P\big\{ |X(w_i)-X(x_i)|  \leq 2r \; \big| \; X(x_1),\ldots,X(x_{i-1}), X(y)\big\}
&\leq c_1\, r^s \psi(|w_i-x_i|)^{-s/2}\\
&\le c_1 r^s |w_i-x_i|^{-\alpha s}
\end{split}
\end{equation}
for each $i=2,\ldots,n$  and $0<r \leq r_0$, where the last inequality follows
from the fact that $\psi(r) \ge r^{2 \alpha}$.

Starting with $\P\big\{|X(y)-X(x_1)|  \leq 2r\big\} \leq c_1 r^s |y-x_1|^{-\alpha s}$,
and applying the conditional probabilities in (\ref{Eq:conp}) inductively, we obtain
\begin{equation*}
\begin{split}
\P\big\{
 | &X(y)-X(x_1)|  \leq r, |X(y)-X(x_2)|  \leq r,\ldots,|X(y)-X(x_n)| \leq r\big\}\\
&\leq
\P\big\{ |X(y)-X(x_1)|  \leq 2r, |X(w_2)-X(x_2)|  \leq 2r,\ldots, |X(w_n)-X(x_n)| \leq 2r\big\}\\
& \leq  (c_1)^n  r^{ns} |y-x_1|^{-\alpha s} |w_2-x_2|^{-\alpha s}\ldots |w_n-x_n|^{-\alpha s}\\
& \leq  (c_1)^n (c_0)^{n \alpha s}   r^{ns}d_a(y,x_1)^{-\alpha s} d_a(w_2,x_2)^{-\alpha s}
\ldots  d_a(w_n,x_n)^{-\alpha s}\\
& =  c_2\,r^{ns}\phi_a(x_1,x_2,\ldots,x_n,y)^{\alpha s},
\end{split}
\end{equation*}
using the definitions   (\ref{transum}) and  (\ref{multpotdef}) of $d_a$ and $\phi_a$
and the choice of the $w_i$.

Inequality (\ref{probdists5}) is immediate from (\ref{probdists4}).
$\Box$

\begin{prop}\label{propexp}
Let $n \geq 1$ and $1 < q \leq n+1$. Then for all $0<s \leq d$,
there exist numbers $c_3>0$ and $r_0>0$
such that for all $0<r\leq r_0$,
\begin{equation} \label{probest2}
\begin{split}
&\E \int \mu_X (B(z,r))^{q-1}  d \mu_X (z) \\
&\quad \leq c_3 r^{s(q-1)} \sum_{a{ \in \cal A}_m} \int\bigg[
\int \cdots \int
\phi_a(x_1,\ldots,x_n,y)^{\alpha s} d\mu(x_1)\ldots d\mu(x_n)\bigg]^{(q-1)/n} d\mu(y).
\end{split}
\end{equation}
\end{prop}

\noindent{\it Proof.}
First note that, for every $y \in {\mathbb R}^N$, using Fubini's theorem and
(\ref{probdists5}), we obtain
\begin{equation} \label{emuomega}
\begin{split}
\E\big(  \mu_X &(B(X(y), r))^n\big)   \\
& = \E\big(\mu\{x :  |X(y)-X(x) | \leq r\}^n \big)   \\
&= \int\cdots\int \P\Big\{ |X(y)-X(x_1)|  \leq r, |X(y)-X(x_2)|  \leq r,   \\
& \hspace{5cm}\ldots,|X(y)-X(x_n)| \leq r\Big\} d\mu(x_1)\ldots d\mu(x_n) \\
&\leq c_2 r^{ns}\int\cdots\int \sum_{a{ \in \cal A}_m}
\phi_a(x_1,x_2,\ldots,x_n,y)^{\alpha s} d\mu(x_1)\ldots d\mu(x_n).
\end{split}
\end{equation}
Since $n/(q-1) \geq 1$, Jensen's inequality,
(\ref{emuomega}) and the power-sum inequality give
\begin{align*}
\E \int \mu_X & (B(z,r))^{q-1}  d \mu_X (z) \\
& = \E \int  \mu_X  (B(X(y),r))^{q-1}  d\mu(y) \\
& \leq  \int\Big[\E\big( \mu_X   (B(X(y),r))^{n}\big)  \Big]^{(q-1)/n}d\mu(y) \\
&\leq c_3 r^{s(q-1)}
 \int\bigg[ \sum_{a{ \in \cal A}_m}\int\cdots\int\phi_a(x_1,x_2,\ldots,x_n,y)^{\alpha s}
 d\mu(x_1)\ldots d\mu(x_n) \bigg]^{(q-1)/n} d\mu(y)\\
 &\leq c_3 r^{s(q-1)}
\int \sum_{a{ \in \cal A}_m} \bigg[ \int\cdots\int\phi_a(x_1,x_2,\ldots,x_n,y)^{\alpha s}
d\mu(x_1)\ldots d\mu(x_n) \bigg]^{(q-1)/n} d\mu(y),
 \end{align*}
to give (\ref{probest2}).
 $\Box$

We now derive the almost sure lower bound for $\lqd (\mu_X)$.

\begin{prop}\label{corlbd}
Let $X: \R^N \to  \R^d $ be the Gaussian random field defined by $(\ref{def:X})$ and
assume that the associated random field $X_0$ satisfies Condition {\bf (C2)}.
Let $\mu$ be a Borel probability measure on $\bbbr^N$ with compact support.
Then for all $q>1$,
\begin{equation}\label{lbdas}
\lqd (\mu_X)   \geq \min\Big\{d, \, \frac 1{\alpha^*}\lqd (\mu) \Big\} \quad \mbox{a.s.}
\end{equation}
\end{prop}

\noindent{\it Proof.}
Without loss of generality we may assume that the support of $\mu$ lies in the
cube $[0,\frac{1}{2})^N$. As before, for each $a \in {\cal A}_m$ we write
${\cal C}_{k}^a$ for the $k$th level cubes in the hierarchy of $m$-ary cubes that
define the metric $d_a$.

Let $\alpha> \alpha^*$ and $0<s <  \min\{d, \lqd (\mu) /\alpha\}$. From (\ref{3.b})
$$\liminf_{k \to \infty} \frac{\log \sum_{C\in {\cal C}_{k}^a }\mu (C)^q}{(q-1)\log m^{-k}}
=\lqd (\mu) >s\alpha$$
so
$$
\limsup_{k \to \infty} \frac{\log \big[m^{s(q-1)\alpha k} \sum_{C\in {\cal C}_{k}^a }
\mu (C)^q\big]}{k} <0.
$$
Note that this estimate holds for all $a \in {\cal A}_m$ since the definition of the
lower generalized dimension (\ref{3.b}) is independent of the origin selected for the
mesh cubes used for the moment sums.

Let $n$ be the integer such that $n \leq q < n+1$. For each $a \in {\cal A}_m$ in turn,
identify the cubes of ${\cal C}_{k}^a$ with the $k$th level vertices of the $M$-ary tree
of Section 3 in the natural way, where $M=m^N$. Thus, with $\bi_j\in I_\infty$ identified
with $\bi_j\in \bbbr^N$ and $\bj\in I_\infty$  with
$\bj\in \bbbr^N$, we have that
$l_i = k(C_i),\, i=1,\ldots, n$, are the levels both of the cubes
and the equivalent vertices in the tree $I_\infty$  in the join set
of $\bi_1,\ldots,\bi_n,\bj$. Setting
$f(l) = m^{\alpha s l}$ in (\ref{defphi}) and
using (\ref{multpotdef}), we get
$$
\phi(\bi_1,\ldots,\bi_n,\bj) =m^{\alpha s  l_1} \cdots  m^{\alpha s l_n} =
m^{\alpha s k(C_1)} \cdots  m^{\alpha s k(C_n)} =
\phi_a(\bi_1,\ldots,\bi_n,\bj)^{\alpha s}.$$
 Thus, Proposition \ref{propexp}  together with Theorem \ref{mainint}  gives that
\be
\E \int   \mu_X (B(z,r))^{q-1}  d\mu_X(z)  \leq c_4 r^{s(q-1)} \label{expmom}
\ee
for all  $r\leq 1$, for some constant $0 <c_4< \infty$.

For all $0<t<s< \min\{d, \lqd (\mu) /\alpha\}$, summing (\ref{expmom}) over
$r= 2^{-k}, k=0,1,2,\ldots$, gives
$$
\E\bigg( \sum_{k=0}^\infty 2^{kt(q-1)}\int \mu_X  (B(z,2^{-k}))^{q-1} d\mu_X (z)\bigg)
\leq  c_4 \sum_{k=0}^\infty 2^{-k(s-t)(q-1)} < \infty.
$$
Thus the bracketed series on the left converges almost surely, so as the generalized
dimensions (\ref{deflgd}) are determined by the sequence of scales $r= 2^{-k}, k=0,1,2,\ldots$,
we conclude that $\lqd (  \mu_X)   \geq  t $ almost surely  for all $0<t < \min\{d, \lqd
(\mu)/\alpha\}$. Since $\alpha > \alpha^*$ is arbitrary, (\ref{lbdas}) follows.
$\Box$

\medskip

\noindent{\it Proof of Theorem \ref{main}.}
The upper bound in (\ref{Eq:main1}) follows from Proposition \ref{Prop:UP},
and the lower bound in (\ref{Eq:main1}) follows from Proposition \ref{corlbd}. This
proves (i).
Part (ii) also follows from Propositions \ref{Prop:UP} and \ref{corlbd}.
Finally, (iii) follows from (ii).
$\Box$

\section{Further remarks}
\setcounter{equation}{0}
\setcounter{theo}{0}
Here are some open problems and remarks about generalized dimensions of random fields
which are not covered by this paper.
\begin{itemize}
\item[1.] Note that (\ref{Eq:main2}) in Theorem \ref{main}  only provides an upper
bound for $\uqd (\mu_X)$. While we believe that if the generalized $q$-dimension of
$D_q (\mu)$ exists then
\[
\uqd (\mu_X)
 = \min\Big\{d, \, \frac 1{\alpha_*}D_q(\mu) \Big\}
\quad \hbox{a.s.},
\]
we have not been able to prove it, because the last inequality in (\ref{Eq:conp}) fails
when $\alpha_* <\alpha < \alpha^*$.

\item[2.] Besides fractional Brownian motion, another important Gaussian random
field is the Brownian sheet $W: \bbbr^N_+ \to \bbbr^d $, $W(x) = (W_1(x),\ldots,W_d(x))$,
which is a centered Gaussian random field with  covariance function given by
$$
\E\big[W_i(x)W_j(y)\big] =  \delta_{ij} \prod_{k=1}^N x_k \wedge
y_k, \quad \  x= (x_1,\ldots,x_N), y= (y_1,\ldots,y_N) \in \bbbr^N_+,
$$
where $\delta_{ij} = 1$ if $i =j$ and $0$ if $i \ne j$, see \cite{Kho}. The
Brownian sheet $W$ is not strongly locally nondeterministic, but satisfies a weaker
form of local nondeterminism, namely, sectorial local nondeterminism as it is called in
\cite{KhoX}. We expect that the conclusion of Corollary \ref{Coro:FBM} still holds
for the Brownian sheet, but since strong local nondeterminism played an important
r\^{o}le in Section 5, a different method may be needed to study the effect of the
Brownian sheet on generalized dimensions.

\item[3.] In recent years several authors have constructed and investigated anisotropic
random fields, see \cite{Xi09} and references therein. Random fractal images under
anisotropic random fields have a  richer geometry than isotropic random fields such
as fractional Brownian motion. It is not clear to what extent our arguments can be
modified for non-isotropic Gaussian fields. Furthermore, it is not clear to what
extent our approach can be used for non-Gaussian random fields,  such as linear
and harmonizable fractional stable random fields, see \cite{Xi11}. For example, difficulties
arise from possible discontinuities of the sample paths.

\item[4.] The case of small moments, that is for $0<q < 1$, is interesting but
is likely to need very different methods.  The negative exponent that appears in 
the generalized dimension integral $\int
\mu(B(x,r))^{q-1}d\mu(x)$  makes such integrals difficult to estimate, 
for example H\"{o}lder's inequality, used throughout the proof of Proposition \ref{fracest}, cannot be used with negative powers.
Extrapolating from the case $q=0$ (when
the generalized dimension reduces to box counting dimension of the support of the
measure), it is possible that some kind of generalized dimension profile is needed,
see \cite{FH,Xi1}.

\end{itemize}

\end{document}